\documentclass{amsart}

\usepackage{latexsym, amsthm}

\usepackage{amsfonts, amsmath, amssymb}

\usepackage{hyperref}

\usepackage{euscript,mathrsfs}

\newtheorem{theorem}{Theorem}[section]
\newtheorem{lemma}[theorem]{Lemma}
\newtheorem{proposition}[theorem]{Proposition}
\newtheorem{conjecture}[theorem]{Conjecture}
\newtheorem{corollary}[theorem]{Corollary}

\theoremstyle{definition}
\newtheorem{definition}[theorem]{Definition}
\newtheorem{example}[theorem]{Example}

\theoremstyle{remark}
\newtheorem{remark}[theorem]{Remark}

\def\Mon{\mathbb M}
\def\Monb{\overline{\mathbb M}}
\def\Nom{\mathbb H}
\def\Fq{{\mathbb F}_q}
\def\s{\mathcal S}
\def\t{\mathcal T}
\def\u{\mathcal U}
\def\m{\mathcal M}
\def\L{\mathcal L}
\def\monord{\succ}
\def\lm{\mathrm{lm}}
\def\x{\mathcal X}
\def\p{\mathsf{p}}
\def\PP{{\mathbb P}}
\def\f{\sigma}
\def\FP{\mathsf {FP}}
\def\SH{\mathsf {SH}}

\def\ZZ{\mathbb{Z}}
\def\AA{{\mathbb A}}
\def\PP{{\mathbb P}}
\def\RM{\mathrm{RM}}
\def\PRM{\mathrm{PRM}}
\def\ev{\mathrm{ev}}

\def\supp{\mathrm{Supp}}
\def\lex{\prec_{\mathrm{lex}}}
\def\lexeq{\preccurlyeq_{\mathrm{lex}}}
\def\glex{\succ_{\mathrm{lex}}}
\def\glexeq{\succcurlyeq_{\mathrm{lex}}}

\def\a{{\alpha}}

\def\FF{{\mathbb F}}

\def\hp3{\widehat{\mathbb P}^3}
\newcommand{\w}{\mathrm{w}}

\def\Aff{{\mathbb A}}

\def\ZZ{{\mathbb Z}}

\title[Combinatorial Approach to Equations over Finite Fields]{A Combinatorial Approach to the Number of Solutions of Systems of Homogeneous Polynomial  Equations  over Finite Fields}

\begin{document}

 \author{Peter Beelen}
\address{
Department of Applied Mathematics and Computer Science,\newline \indent
Technical University of Denmark, 
2800 Kgs. Lyngby, Denmark}
\email{pabe@dtu.dk}
\author{Mrinmoy Datta}
\address{
Department of Mathematics,
Indian Institute of Technology Hyderabad,\newline \indent
Kandi, Sangareddy, Telanagana,   502285, India}
\email{mrinmoy.datta@math.iith.ac.in} 

\author{Sudhir R. Ghorpade}
\address{Department of Mathematics, %
Indian Institute of Technology Bombay,\newline \indent
Powai, Mumbai 400076, India}
\email{srg@math.iitb.ac.in}


\date{}

\begin{abstract}
We give a complete conjectural formula for the  number $e_r(d,m)$ of maximum possible $\Fq$-rational points on a projective algebraic variety defined by 
$r$ linearly independent homogeneous polynomial equations of degree $d$ in $m+1$ variables with coefficients in the finite field $\Fq$ with $q$ elements, when $d<q$. It is shown that this formula holds in the affirmative for several values of $r$. In the general case, 
we give explicit lower and upper bounds for $e_r(d,m)$ and show that they are sometimes attained.
Our approach uses a relatively recent result, 
called the projective footprint bound, together with results from extremal combinatorics such as the Clements-Lindstr\"om Theorem and its variants. Applications to the problem of determining the generalized Hamming weights of projective Reed-Muller codes are also included.
\end{abstract}

\maketitle

\section{Introduction}
\label{sec:intro}
Fix a prime power $q$ and positive integers $r,d,m$. Let $\Fq$ denote the finite field with $q$ elements and $\Fq[x_0, \dots , x_m]$ the polynomial ring in $m+1$ variables $x_0, x_1, \dots , x_m$ with coefficients in $\Fq$. For any homogeneous polynomials $F_1, \dots , F_r$ in $\Fq[x_0, \dots , x_m]$, let
$V(F_1, \dots , F_r)$ denote the closed subvariety of the projective $m$-space $\PP^m$ (over an algebraic closure of $\Fq$) given by the vanishing of $F_1, \dots , F_r$, and let $V(F_1, \dots , F_r)(\Fq)$ be the set of its $\Fq$-rational points,
i.e., the set of all $\Fq$-rational common zeros in $\PP^m$ of $F_1, \dots , F_r$. Define
\begin{equation}
\label{erdm}
e_r(d,m) := \max_{F_1, \dots , F_r} \left| V(F_1, \dots , F_r)(\Fq) \right|,
\end{equation}
where the maximum is over all possible families $\{F_1, \dots , F_r\}$ of $r$ linearly independent homogeneous polynomials of degree $d$ in $ \Fq[x_0, \dots , x_m]$. Note that the condition on linear independence implies that $r$ can be at most $\binom{m+d}{d}$. Note also that an obvious upper bound for $e_r(d,m)$ is $\p_m$,  where for $j\in \ZZ$, by $\p_j$ we denote $|\PP^j(\Fq)|$, i.e., $\p_j:= q^j + q^{j-1} + \dots + q + 1$ if $j\ge 0$, and 
$\p_j:=0$ if $j< 0$.

Explicit determination of $e_r(d,m)$ is an open problem, in general,  and it has been of some interest for about two decades. 
While it is easy to see that $e_r(1,m) = \p_{m-r}$ for $r\le m+1$ and $e_r(d,1)= d-r+1$ for $r\le d+1 \le q$ (see, e.g., \cite[\S\,2.1]{DG2}), the case of $r=1$ is rather nontrivial. Here it was conjectured by M. Tsfasman that
\begin{equation}
\label{TConj}
e_1(d,m) = dq^{m-1} + \p_{m-2}  \quad \text{for } d \le q.
\end{equation}
This was subsequently proved by Serre \cite{Se} and, independently, by S{\o}rensen \cite{So} in 1991.
In the general case, an intricate formula for $e_r(d,m)$ for $d< q-1$ was conjectured by Tsfasman and Boguslavsky (cf. \cite{Bog, DG}), and this was proved in the affirmative by Boguslavsky \cite{Bog} in 1997 for $r=2$. The case of $r>2$ remained open for a considerable time.
Eventually, it was proved in \cite{DG2} and \cite{DG} that the conjectural formula of Tsfasman and Boguslavsky is true if $r\le m+1$, and it can be false if $r> m+1$. In \cite{DG2}, a new conjectural formula for $e_r(d,m)$ was proposed for many (but not all) values of $r$, namely  
for $r\le \binom{m+d-1}{d-1}$. We will refer to this as the \emph{incomplete conjecture}. To state it, let us first define an important combinatorial quantity whose genesis lies in the work of Heijnen and Pellikaan \cite{HP} related to an affine counterpart of the problem of finding $e_r(d,m)$. 
For $1\le r \le \binom{m+d}{d}$, we define
$$
H_r(d,m):= \sum_{i=1}^m \alpha_i q^{m-i}, \quad \text{where $(\alpha_1, \dots, \alpha_m)$ is the $r^{\rm th}$ element of } Q^m_{\le d},
$$
and where $Q^m_{\le d}$ denotes the collection, ordered in descending lexicographic order, of all
$m$-tuples $(\beta_1, \dots , \beta_m)$ of integers satisfying $0\le \beta_i <q$ for $i=1, \dots , m$ and $\beta_1+\cdots + \beta_m \le d$. For example, if $d<q$, then $H_1(d,m) = dq^{m-1}$ and $H_{\binom{m+d}{d}}(d,m)=0$,
since $(d,0, \dots , 0)$ and $(0,0, \dots , 0)$ are clearly the first and the last $m$-tuples of $Q^m_{\le d}$ 
in descending lexicographic order. As a convention, we set
\begin{equation}
\label{convention}
H_0(d,m) := q^m \text{ for } d,m\ge 0
 \quad \text{and} \quad H_1(d,m) := 0  \text{ if $d=0$ or $m=0$}.
\end{equation}
In this way, $H_r(d,m)$ is defined for all nonnegative integers $r,d,m$ with $r \le \binom{m+d}{d}$.
The ``incomplete conjecture" of 
\cite{DG2} can now be stated as follows.
\begin{equation}
\label{IC}
e_r(d,m) = H_r(d-1,m) + \p_{m-1} \quad \text{ for $1\le r \le\binom{m+d-1}{d-1}$ and $1\le d < q$.}
\end{equation}
For example, if $r=1$, then this says that $e_1(d,m) = (d-1)q^{m-1} + \p_{m-1}$, which agrees with the Serre-S{\o}rensen formula \eqref{TConj}. Note also that \eqref{IC} holds trivially when $d=1$ or $m=1$. Results of \cite{DG2} prove \eqref{IC} in the affirmative if $r\le m+1$ and $d < q-1$. The validity of \eqref{IC} was extended further in \cite{BDG} to $r\le \binom{m+2}{2}$ and $1< d < q$. This, then, is currently the best known general result as far as an explicit determination of $e_r(d,m)$ is concerned. Apart from this, the last few values of $e_r(d,m)$ were determined in \cite{DG3} using the connection with coding theory (explained in Section \ref{sec:codes})
 and the work of S{\o}rensen \cite{So};
more precisely, it is shown in \cite[Thm. 4.7]{DG3} that
\begin{equation}
\label{LastFew}
e_{\binom{m+d}{d} - t} (d,m) = t \quad \text{for } t=0,1, \dots , d.
\end{equation}

We are now ready to describe the main results of this paper. First, we extend \eqref{IC} to a conjectural formula for $e_r(d,m)$ for all permissible values of $r,d,m$ with $d<q$. To state this ``complete conjecture'', let us first observe that
\begin{equation}
\label{Pascal}
\binom{m+d}{d} = \binom{m+d-1}{d-1} + \binom{m+d-2}{d-1} + \cdots + \binom{d-1}{d-1}
\end{equation}
and that 
for any positive integer $r <  \binom{m+d}{d}$, there are unique integers $i,j$ 
such that
$$
r = \binom{m+d-1}{d-1} 
+ \cdots + \binom{m+d-i}{d-1} + j, 
\ \,
0\le i \le m, \text{ and } 0 \le j <  \binom{m+d-i-1}{d-1}.
$$
By convention, and in accordance with \eqref{Pascal}, we set $i:=m$ and \hbox{$j:=  \binom{m+d-i-1}{d-1}=1$} when $r =  \binom{m+d}{d}$. With $i$ and $j$ thus defined (for a given value of $r$),
the ``complete conjecture" states that
\begin{equation}
\label{CC}
e_r(d,m) = H_j(d-1,m-i) + \p_{m-i-1} \quad \text{ for $1\le r \le\binom{m+d}{d}$ and $1\le d < q$.}
\end{equation}
Note that if $r < \binom{m+d-1}{d-1}$, then $i=0$ and $j=r$, whereas if $r=\binom{m+d-1}{d-1}$, then
$i=1$ and $j=0$. Thus \eqref{CC} reduces to \eqref{IC} in this case, thanks to our conventions.  In particular, from \cite[Thm. 5.3]{BDG}, we see that \eqref{CC} holds in the affirmative if $r\le \binom{m+2}{2}$.
We provide further evidence for the ``complete conjecture" in this paper by showing that it holds in the affirmative for
an additional $md$ values of $r$, namely for
$$
r =\! \binom{m+d-1}{d-1}  + \cdots + \binom{m+d-i}{d-1} - t \quad \text{where } 1\le i \le m \text{ and } 0 \le t \le d-1;
$$
in fact, for $r$ as above, we obtain $e_r(d,m) =  \p_{m-i} + t$. These results 
are also valid when $i=m+1$, but in view of \eqref{Pascal}, this case is 
already  covered by  \eqref{LastFew}.	
In the general case, we show that the conjectural formula is always a lower bound even when $d=q$, that is,
\begin{equation}
 \label{LB}
e_r(d,m) \ge H_j(d-1,m-i) + \p_{m-i-1} \quad \text{for $1\le r \le\binom{m+d}{d}$ and $1\le d \le q$.}
\end{equation}
The conjectural formula \eqref{CC} as well as the 
lower bound \eqref{LB} for $e_r(d,m)$ can be  described by the alternative formula 
$$
\p_{s_d-d} +  \lfloor q^{s_{d-1}-d+1} \rfloor+ \lfloor q^{s_{d-2}-d+2} \rfloor + \cdots + \lfloor q^{s_1-1} \rfloor, 
$$
where $s_1, \dots , s_d$ are unique integers satisfying the $d$-binomial expansion
$$
\binom{m+d}{d} - r = \binom{s_d}{d} + \binom{s_{d-1}}{d-1} + \cdots + \binom{s_1}{1} \quad \text{and} \quad 
s_d > s_{d-1} > \cdots > s_1 \ge 0.
$$

We also find in this paper an explicit upper bound for $e_r(d,m)$ using methods from extremal combinatorics and a projective counterpart of $H_r(d,m)$, that we denote by $K_r(d,m)$. More precisely, we show that 
for $1\le r \le\binom{m+d}{d}$ and $1\le d < q$,
\begin{equation}
 \label{UB}
e_r(d,m) \le K_r(d,m), \quad \text{where} \quad K_r(d,m) := \sum_{i=0}^{m}  a_i \p_{m-i-1},
\end{equation}
and where $(a_0, a_1, \dots , a_m)$ is the $r$-th element, in descending lexicographic order, of the set of all $(m+1)$-tuples $(b_0, b_1, \dots , b_m)$ of nonnegative integers satisfying $b_0+b_1+ \cdots + b_m = d$. It is also shown that this upper bound is attained for several values of $r$. In turn, this plays a crucial role in ascertaining the validity of \eqref{CC} for the additional $md$ values of $r$ mentioned earlier.

The determination of $e_r(d,m)$ is directly related to the determination of the generalized Hamming weights (also known as higher weights) of 
the projective Reed-Muller codes $\PRM_q(d,m)$, which go back to Lachaud \cite{La}. This connection has been explained in \cite[\S\,4]{DG3} when $d\le q$. We elucidate it further in Section \ref{sec:codes} by noting that in general (when $d$ can be larger than $q$), it is more natural to consider a variant of $e_r(d,m)$, called $\overline{e}_r(d,m)$, wherein one takes into account the vanishing ideal of $\PP^m(\Fq)$. This also brings to the fore a basic notion of projective reduction that was enunciated in \cite{BDG18}. We remark that the affine counterpart of the problem of determining $\overline{e}_r(d,m)$ corresponds to determining the generalized Hamming weights of Reed-Muller codes $\RM_q(d,m)$, and this has been solved by V. Wei \cite{W} when $q=2$, and by Heijnen and Pellikaan \cite{HP}, in general. (See also \cite{BD}.)

The methods used in proving the main results of this paper differ significantly from those in our earlier works such as \cite{DG2} and \cite{BDG}. Here we adopt a combinatorial approach and the groundwork for this has been laid in \cite{BDG18} where a \emph{projective footprint bound} for the number of $\Fq$-rational points of arbitrary projective algebraic varieties defined over $\Fq$ was obtained. This groundwork is combined in this paper with methods from extremal combinatorics and a culmination of ideas such as the Kruskal-Katona Theorem, a lemma of Wei \cite[Lem. 6]{W}, the Clements-Lindstr\"om Theorem, and a theorem of Heijnen \cite[Appendix A]{H} (see also \cite[Thm. 5.7]{HP} and \cite[Thm. 3.8]{BD}).
These feed into the results in Sections \ref{combin}, \ref{combin2} and \ref{sec:er} that form the technical core of this paper. 
For a leisurely introduction to extremal combinatorics and some of the classical results mentioned above, one  may refer to
the book of Anderson \cite{A}. 

%

\section{A Lower Bound and a Conjecture}
\label{sec:new}

In this section, we shall prove the lower bound \eqref{LB} and formally state our ``complete conjecture" concerning $e_r(d,m)$. Recall that $e_r(d,m)$  is defined by \eqref{erdm} for positive numbers $r,d,m$ with $r \le \binom{m+d}{d}$. One can extend the definition to include the trivial cases when $r$ or $d$ or $m$ is zero, or simply, set the following convention.
 \begin{equation}
\label{convention-erdm}
e_0(d,m) := \p_m \text{ for } d,m\ge 0
 \quad \text{and} \quad e_1(d,m) := 0  \text{ if $d=0$ or $m=0$}.
\end{equation}
Thus $e_r(d,m)$ is defined for all nonnegative integers $r,d,m$ satisfying $r \le \binom{m+d}{d}$.
%

In the remainder of this section, $m$ denotes a positive integer and, as usual, $\Fq$ denotes the finite field with $q$ elements.

\subsection{Zeros of Affine Varieties over Finite Fields}\label{subsec:Hr1}
Let $r, d$ be any nonnegative integers. 
Given any polynomials $f_1, \dots , f_r\in \Fq[x_1, \dots , x_m]$, we shall denote
by $Z(f_1, \dots, f_r)(\Fq)$ the set of all $\Fq$-rational points of the affine algebraic variety in $\Aff^m$ (over an algebraic closure of $\Fq$) defined by $f_1, \dots , f_r$; in other words,
$$
Z(f_1, \dots, f_r)(\Fq) =\{ (a_1,\dots , a_m)\in \Fq^m : f_j(a_1,\dots , a_m) =0 \text{ for all } j=1,\dots ,r\}.  
$$
Note that if $r=0$, then $Z(f_1, \dots, f_r)(\Fq) =\Fq^m$. 
We shall now define a natural affine analogue of  $e_{r}(d,m)$.  
For $0\le r \le \binom{m+d}{d}$, we define 
$$
e^{\Aff}_r(d,m):=\max_{f_1, \dots, f_r} |Z(f_1, \dots, f_r)(\Fq)|,
$$
where the maximum is taken over families of $r$ linearly independent polynomials $f_1, \dots , f_r$ of degree $\le d$ in $\Fq[x_1, \dots , x_m]$.

As explained in \cite[\S 2.1]{BDG}, the result of Heijnen and Pellikaan \cite[Thm. 5.10]{HP} in the case $d<q$ can be stated as follows. Here $H_r(d,m)$ is as defined in the introduction, including the conventions given in \eqref{convention}.

\begin{theorem}[Heijnen-Pellikaan]
\label{Hr1}
If $0\le d < q$ and $0\le r \le
\binom{m+d}{d} $, then
$$
{e}^{\Aff}_r(d,m)  = H_r (d, m).
$$
\end{theorem}

A more general version of this result will be discussed later (in \S\,\ref{subsec:RedEqns}).

\subsection{A Lower Bound for $e_r(d,m)$}
We begin by noting a simple and well-known fact whose proof is outlined for the sake of completeness.

\begin{lemma}\label{lem:binom}
Let $d$ be a positive integer. Then
\begin{equation}
\label{Pascal2}
\binom{m+d}{d} = \sum_{a=1}^{m+1}\binom{m+d-a}{d-1} .
\end{equation}
Moreover, for any nonnegative integer $r <  \binom{m+d}{d}$, there are unique integers $i,j$ with
\begin{equation}
\label{BinomExp}
r = j + \sum_{a=1}^{i}\binom{m+d-a}{d-1} , \quad
0\le i \le m, \quad \text{and} \quad 0 \le j <  \binom{m+d-i-1}{d-1}.
\end{equation}
\end{lemma}

\begin{proof}
The identity in \eqref{Pascal2} follows easily from induction on $m$. 
If $0\le r <  \binom{m+d}{d}$, then  the largest nonnegative integer $i$ such that $\sum_{a=1}^{i}\binom{m+d-a}{d-1} \le r$ clearly satisfies $0\le i \le m$, thanks to \eqref{Pascal2}.
Thus
\eqref{BinomExp} holds with $j: = r - \sum_{a=1}^{i}\binom{m+d-a}{d-1}$.
\end{proof}

We remark that although $m$ is assumed to be a fixed positive integer, the identity in \eqref{Pascal2} holds trivially also when $m=0$, and this fact may be tacitly assumed in the sequel. Our next result is 
a general lower bound for $e_r(d,m)$ when $d\le q$.

\begin{theorem}\label{lem:lower}
Let $d,r$ be positive integers with
$d \le q$ and $r \le \binom{m+d}{d}$,
and let $i,j$ be as in 
\eqref{BinomExp}
if 
$r < \binom{m+d}{d}$, while  $i:=m$ and $j:=  \binom{m+d-i-1}{d-1}=1$ if 
$r = \binom{m+d}{d}$. Then
$$
e_r(d,m) \ge H_j(d-1,m-i)+\p_{m-i-1}.
$$
\end{theorem}

\begin{proof}
If $r = \binom{m+d}{d}$, then clearly, $e_r(d,m) =0$ and 
$H_1(d-1, m-m) + \p_{m-m-1} = 0$, 
as per the conventions in equation \eqref{convention}. Now assume that
$r < \binom{m+d}{d}$, and let $i, j$ be as in \eqref{BinomExp}. We shall prove the desired inequality
by producing a set $B$ 
of $r$ linearly independent polynomials in $\Fq[x_0, \dots, x_{m}]_{d}$ with the property 
that $|V(B)(\Fq)| = H_j(d-1,m-i)+\p_{m-i-1}$.

First, for each positive integer
$a\le i$, let $\mathcal{B}_{a}$ be a basis of the $\Fq$-vector space $x_{m-a+1}\Fq[x_0,\dots,x_{m-a+1}]_{d-1}$; for instance, 
$\mathcal{B}_{a}$ can be the set of monomials of
degree $d$ in $x_0,\dots,x_{m-a+1}$ that are divisible by $x_{m-a+1}$. Clearly, the sets $\mathcal{B}_{a}$ are disjoint and $|\mathcal{B}_{a}| = \binom{m+d-a}{d-1}$ for $1\le a \le i$. Note ~that
$$
V\Big(\bigcup_{a=1}^i \mathcal{B}_{a} \Big)(\Fq)=\{(a_0:\cdots:a_m) \in \PP^m(\Fq) \, | \, a_m=\cdots=a_{m-i+1}=0\},
$$

Next,  since $d-1< q$, by Theorem \ref{Hr1}, we obtain $j$ linearly independent polynomials $f_1, \dots, f_j$ each
of degree 
at most $d-1$
in the polynomial ring $\Fq[x_0, \dots, x_{m - i - 1}]$ such that
$
|Z_{m-i}(f_1, \dots, f_j) (\Fq)| = H_j (d-1, m - i),
$
where $Z_{m-i}(f_1, \dots, f_j)$ denotes the set of common zeroes of $f_1, \dots, f_j$ in the $(m-i)$ dimensional affine subspace  of $\PP^m$ given by  $\{ [a_0:\cdots:a_m] : a_{m-i} =1, \; a_s = 0 \ \text{for }
m-i< s \le m\}$. 
Let $F_1, \dots , F_j$ be the 
polynomials 
obtained, respectively, by homogenizing $f_1, \dots , f_j$ to degree $d$ with respect to the variable $x_{m-i}$.
Clearly, $F_1, \dots , F_j \in \Fq[x_0, \dots, x_{m - i }]_d$~and  they are linearly independent.
Also, $\Fq[x_0, \dots, x_{m - i }] \cap \mathcal{B}_{a}$ is empty
for each $a=1, \dots , i$. Consequently, $B$ is a linearly independent subset of $\Fq[x_0, \dots, x_{m }]_d$~and
$$
|B| =  j + \sum_{a=1}^i \binom{m+d-a}{d-1}  = r \quad \text{where} \quad
B:= \{F_1, \dots , F_j\} \cup \Big(\bigcup_{a=1}^i \mathcal{B}_{a} \Big),
$$
with the 
convention that the relevant sets are empty if $j=0$ or $i=0$.
Further, by intersecting $V(B) = V(F_1, \dots , F_j) \cap V(\cup_{a=1}^i \mathcal{B}_{a})$ with the affine patch $\{x_{m-i}=1\}$ and the hyperplane $\{x_{m-i}=0\}$ of $\PP^m$, we see that
$|V(B)(\Fq)|$ equals
$$
|Z_{m-i}(f_1, \dots, f_j) (\Fq)|  + \left| \{ [a_0:\cdots:a_m] \in \PP^m(\Fq) : a_m=\cdots=a_{m-i}=0\} \right|,
$$
which is $H_j(d-1,m-i)+\p_{m-i-1}$, as desired.
\end{proof}

We remark that when $d=q$, the lower bound in Theorem \ref{lem:lower} is not attained, in general. This is shown in \cite[\S\,6]{BDG}, where exact values of $e_r(q,m)$ are obtained for $1\le r \le m+1$. 

\subsection{The ``Complete Conjecture"}
It appears plausible that the lower bound in Theorem \ref{lem:lower} is attained when $d< q$. More precisely, we conjecture the following.

\begin{conjecture}\label{conj:extended}
Let $d,r$ be 
integers with
$1\le d < q$ and $1\le r \le \binom{m+d}{d}$,
and let $i,j$ be as in 
\eqref{BinomExp}
if 
$r < \binom{m+d}{d}$, while  $i:=m$ and $j:=  \binom{m+d-i-1}{d-1}$ if 
$r = \binom{m+d}{d}$. Then
$$
e_r(d,m)=H_j(d-1,m-i)+\p_{m-i-1}.
$$
\end{conjecture}

As noted in the introduction, if $1 \le r \le \binom{m+d-1}{d-1}$, then Conjecture \ref{conj:extended} reduces to \cite[Conjecture 6.6]{DG2}. The case $d=1$ of Conjecture \ref{conj:extended} holds trivially, while the case $d=2$ follows from the work of Zanella \cite{Z}. Further, it was shown in \cite{BDG} that Conjecture \ref{conj:extended} holds when $2 \le d<q$ and $r \le \binom{m+2}{2}$.
%
%
%

We will now describe an alternative formulation of Conjecture \ref{conj:extended} 
deduced from the alternative description of $H_r(d,m)$ given in \cite{B}. Before stating it, let us recall that given any positive integer $d$, we can express every nonnegative integer $N$ as 
$$
N =\sum_{a=1}^d \binom{s_a}{a}
 \quad \text{for unique $s_a\in \ZZ$ 
with} \quad
s_d > s_{d-1} > \cdots > s_1 \ge 0.
$$
This is called the \emph{$d$-binomial representation} of $N$ or the \emph{$d$-th Macaulay representation} of $N$. 
We will find it convenient to consider $m_a:= s_a -a$ for $1\le a \le d$ so as to write the above expansion for $N$ as 
\begin{equation}\label{eq:MacaulayRep}
N = \sum_{a=1}^d \binom{m_a+a}{a} 
\ \text{ for unique $m_a\in \ZZ$ 
with } \
m_d \ge m_{d-1} \ge \cdots \ge m_1 \ge -1.
\end{equation}
We may refer to $(m_d, \dots , m_1)$ as the \emph{Macaulay $d$-tuple} corresponding to $N$. Observe that if $M$ is any nonnegative integer, then 
\begin{equation}\label{eq:MacaulayCoeff}
0\le N <\binom{M+d}{d} \Longrightarrow 
M -1 \ge m_d \ge m_{d-1} \ge \cdots \ge m_1 \ge -1.
\end{equation}
 The following result is a direct consequence of the proof of \cite[Thm. 3.1]{B} and we remark that its proof
does not use the theorem of Heijnen-Pellikaan (Theorem \ref{Hr}). 

\begin{lemma}\label{lem:pabe}
Assume that $1\le d < q$ and $0 \le r \le \binom{m+d}{d}$. Suppose the Macaulay $d$-tuple corresponding to 
$\binom{m+d}{d} -r$ is $(m_d, \dots , m_1)$. Then 
$$
H_r(d,m) = \sum_{a=1}^d \lfloor q^{m_a} \rfloor.
$$ 
\end{lemma}

\begin{proof}
If $r=0$, then  $H_r(d,m)=q^m$, in accordance with our convention \eqref{convention}. On the other hand, the
Macaulay $d$-tuple corresponding to $\binom{m+d}{d}$ is clearly $(m, -1, \dots , -1)$. So the desired equality holds when $r=0$. 
For $1 \le r \le \binom{m+d}{d}$, the desired equality 
 is a special case of \cite[Thm. 3.1]{B} and its proof, since for $d< q$, the dimension $\rho_q(d,m)$ of the Reed-Muller code $\RM_q(d,m)$ is 
$\binom{m+d}{d}$.
\end{proof}

\begin{corollary}\label{cor:MacExp}
Assume that $1\le d \le q$ and $1\le r \le \binom{m+d}{d}$. Let $i,j$ be as in 
\eqref{BinomExp}
if 
$r < \binom{m+d}{d}$, while  $i:=m$ and $j:=  \binom{m+d-i-1}{d-1}$ if 
$r = \binom{m+d}{d}$. Also let $(m_d, \dots , m_1)$ be the Macaulay $d$-tuple corresponding to 
$\binom{m+d}{d} -r$.  Then
$$
H_j(d-1,m-i)+\p_{m-i-1} = \p_{m_d} + \sum_{a=1}^{d-1} \lfloor q^{m_a} \rfloor.
$$
\end{corollary}

\begin{proof}
First, note that if $d=1$, then $i=r$ and $j=0$, and so in view of \eqref{convention}, $H_j(d-1, m-i) +\p_{m-i-1} =
q^{m-r} + \p_{m-r-1} =\p_{m-r}$. Also,  the Macaulay $1$-tuple corresponding to $m+1-r$  is clearly $(m-r)$. Thus 
the desired equality holds when $d=1$. 
Likewise, if $i=m$, then clearly $j$ is $0$ or $1$ and 
$\binom{m+d}{d} -r$ is $1$ or $0$, thanks to \eqref{Pascal2}; in this case, $m_a=-1$ for $1\le a < d$, whereas $m_d$ is $0$ or $-1$ according as $j$ is $0$ or $1$. Hence, in view of \eqref{convention}, the desired equality also holds when $i=m$. 

Now suppose $d>1$ and $0\le i < m$.    
Then  $1\le d-1 < q$ and $1\le  r < \binom{m+d}{d}$.  
%
Using the equality in \eqref{BinomExp} 
together with \eqref{Pascal2}, we can write 
\begin{eqnarray*}
\binom{m+d}{d} - r & =&  
\binom{m+d-i-1}{d-1} - j 
+ \sum_{a=1}^{m-i} \binom{m+d-i-a-1}{d-1}  \\
 & =&  \binom{m+d-i-1}{d} +  \binom{m+d-i-1}{d-1} - j  ,
\end{eqnarray*}
where the last equality follows from 
 \eqref{Pascal2} with $m$ replaced by $m-i-1$. In case $j=0$, the last expression is simply $\binom{m+d-i}{d}$, and in this case $m_d=m-i$ while $m_a=-1$ for $1\le a < d$. Thus, in view of \eqref{convention}, we see that when $j=0$, 
$$
H_j(d-1,m-i)+\p_{m-i-1} = q^{m-i} +\p_{m-i-1} =\p_{m-i} =
\p_{m_d} + \sum_{a=1}^{d-1} \lfloor q^{m_a} \rfloor,
$$
as desired. Now suppose $0< j < \binom{m+d-i-1}{d-1}$. By \eqref{eq:MacaulayCoeff} and Lemma \ref{lem:pabe}, it follows that if $(\mu_{d-1}, \dots , \mu_1)$ is the Macaulay $(d-1)$-tuple corresponding to $ \binom{m+d-i-1}{d-1} - j$, then 
$m-i-1 \ge \mu_{d-1} \ge \dots \ge \mu_1$, and further, 
\begin{equation}\label{eq:Hj}
\binom{m+d-i-1}{d-1} - j  =  \sum_{a=1}^{d-1} \binom{\mu_a+a}{a} \quad \text{and} \quad 
H_j(d-1,m-i) = \sum_{a=1}^{d-1} \lfloor q^{\mu_a} \rfloor .
\end{equation}
Substituting this in the expression obtained earlier for $\binom{m+d}{d} - r$, we see that 
$$
\binom{m+d}{d} - r  = \binom{m+d-i-1}{d} + \sum_{a=1}^{d-1} \lfloor q^{\mu_a} \rfloor .
$$
This 
together with 
the uniqueness of Macaulay $d$-tuples implies that $m_d = m-i-1$ and $m_a = \mu_a$ for $1\le a < d$. 
Consequently,  
\eqref{eq:Hj} yields the desired equality
\end{proof}

In view of Corollary \ref{cor:MacExp}, 
the lower bound in Theorem \ref{lem:lower} and 
the conjectural formula for $e_r(d,m)$ in Conjecture \ref{conj:extended} can be written as 
$$
\p_{m_d} + \sum_{a=1}^{d-1} \lfloor q^{m_a} \rfloor,
$$
where $(m_d, \dots, m_1)$ is the Macaulay $d$-tuple corresponding to $\binom{m+d}{d} -r$.  
\section{Projective Reduction, Shadows and Footprints}
\label{sec:prelim} 

In this section, we review some preliminary notions and results, which will be useful to us in the remainder of the paper. Throughout this section, $m$ denotes a positive integer and $\Fq$ the finite field with $q$ elements.

\subsection{Projective Reduction}
\label{subsec:projred}
Recall that a monomial $\mu \in \Fq[x_1, \dots, x_m]$ given by $\mu = x_1^{\alpha_1} \cdots x_m^{\alpha_m}$ is said to be \textit{reduced} if $0 \le \alpha_i \le q-1$ for all $i = 1, \dots, m$ and a polynomial $F \in \Fq[x_1, \dots, x_m]$ is said to be reduced if it is an $\Fq$-linear combination of reduced monomials. It is 
well-known (see, e.g., \cite[Ch. 2]{Joly} or \cite{G}) 
that the set of all reduced monomials gives rise to a basis of the $\Fq$-vector space $\Fq[x_1, \dots, x_m] / I (\AA^m (\Fq))$, where $I (\AA^m (\Fq))$ denotes the ideal consisting of all polynomials in $\Fq[x_1, \dots, x_m]$ vanishing at every point of $\AA^m (\Fq)$. More precisely, any element of $\Fq[x_1, \dots, x_m] / I (\AA^m (\Fq))$ can be written uniquely as $\tilde{f} + I (\AA^m (\Fq))$ for some reduced polynomial $\tilde{f}$.

A projective analogue of the above notion and result is given in \cite[\S\,2]{BDG18}.
We recall this below. Here, and hereafter,
we denote by $\Mon$ the set of all monomials in the $m+1$ variables $x_0,\dots, x_{m}.$
\begin{definition}
For a nonnegative integer $a$, 
let $\overline{a}$ be the unique integer satisfying
$$
0 \le \overline{a} < q\quad \text{and} \quad
\overline{a} =
\begin{cases}
0  \ \ \mathrm{if} \  a=0 \\
\tilde{a}  \ \ \mathrm{if} \ a > 0 \ \ \mathrm{and} \  \tilde{a} \equiv a({\rm mod} \ {q-1}), \ \mathrm{where} \ 0 < \tilde{a} \le q-1.
\end{cases}$$
Let $\mu  \in \Mon$. If $\mu \neq 1$, then we may write $\mu =x_0^{a_0}\cdots x_\ell^{a_\ell}$, where $0 \le \ell \le m$ and $a_0, \dots , a_\ell$ are nonnegative integers with $a_\ell>0$. Define
$$
\overline{\mu}:=x_0^{\overline{a_0}} \cdots x_{\ell-1}^{\overline{a_{\ell-1}}} x_\ell^{a_\ell+\sum_{j=0}^{\ell-1}(a_j-\overline{a_j})}.
$$
If $\mu=1$, then we define $\overline{\mu}=1$. Note that $\overline{\mu}\in \Mon$ with $\deg \overline{\mu} = \deg \mu$. We call $\overline{\mu}$ the \emph{projective reduction} of $\mu$.
Any polynomial $F \in \Fq[x_0, \dots, x_m]$ can be written uniquely as $F = \sum_{i=1}^n c_i \mu_i$ where $c_i \in \Fq\setminus\{0\}$ and $\mu_i \in \Mon$. We define $\overline{F}$, the \emph{projective reduction} of $F$, as $\overline{F} = \sum_{i=1}^n c_i \overline{\mu_i}$.
A monomial $\mu \in \Mon$ (resp.~polynomial $F \in \Fq[x_0,\dots,x_m]$) is said to be \textit{projectively reduced} if $\overline{\mu}=\mu$ (resp.~$\overline{F}=F$).
\end{definition}
It is easy to see from the definition that a reduced polynomial is projectively reduced. In particular, a polynomial of degree $d \le q-1$ is necessarily reduced and hence projectively reduced. A polynomial of degree $q$ is not necessarily reduced but is always projectively reduced. Clearly, a polynomial of degree $d >q$ may not even be projectively reduced. It is easy to see that if $F\in \Fq[x_0, \dots , x_m]$, then
\begin{equation}\label{eq:redeval}
F(c_0, \dots , c_m) = \overline{F}(c_0, \dots , c_m) \quad \text{for all } (c_0, \dots , c_m) \in \Fq^{m+1}.
\end{equation}
For further properties of projective reduction and projectively reduced polynomials we refer to \cite[Proposition 2.2]{BDG18}.

Throughout this article, we 
denote by $\Monb$ the set of projectively reduced monomials in $m+1$ variables $x_0, \dots , x_m$. Further, for any nonnegative integer $e$, we denote by $\Monb_e$ the set of all projectively reduced monomials in $\Monb$ of degree $e$. Clearly, $\Monb$ equals the disjoint union $\coprod_{e \ge 0}\Monb_e$.

Let $I(\PP^m (\Fq))$ denote the ideal of $\Fq[x_0, \dots, x_m]$ generated by the 
homogeneous polynomials that vanish at all points of $\PP^m (\Fq)$. It was shown by Mercier and Rolland \cite{MR} that this ideal is equal to  the ideal $\Gamma_q(\Fq)$ of $\Fq[x_0, \dots, x_m]$ generated by
$\{x_i^qx_j-x_ix_j^q \mid 0 \le i < j \le m\}$.
It was further proved in \cite[Theorem 2.8]{BDG18} that the set $\{x_i^qx_j-x_ix_j^q \mid 0 \le i < j \le m\}$ forms a universal Gr\"{o}bner basis for $\Gamma_q(\Fq)$.
To conclude this subsection, we note that the projectively reduced monomials give rise to a  basis of the $\Fq$-vector space $\Fq[x_0, \dots, x_m] / I (\PP^m (\Fq))$.  More precisely, any element of $\Fq[x_0, \dots, x_m] / I (\PP^m (\Fq))$ can be written uniquely as $f + I (\PP^m (\Fq))$, where $f$ is a projectively reduced polynomial.
For a proof, see 
\cite[Corollary 2.10]{BDG18}.

\subsection{Shadows and Footprints}
\label{subsec:pfb}
Given $\s \subseteq \Mon$ and  a nonnegative integer $e$, define
$$
\nabla_e(\s) := \{\mu \in \Monb_e : \nu \mid \mu \ \mathrm{for \ some} \ \nu \in \s\} \ \ \mathrm{and} \ \ \ \Delta_e(\s):= \Monb_e \setminus \nabla_e(\s).
$$
The set $\nabla_e(\s)$ is called the \textit{shadow} of $\s$ in $\Monb_e$, while the set $\Delta_e (\s)$ is known as the \textit{footprint} of $\s$ in $\Monb_e$.

Recall that by a \emph{term order} on $\Mon$, 
one means a total order $\prec$ 
on the set $\Mon$ of all monomials in $x_0,\dots,x_m$ such that (i) $1 \preccurlyeq \mu$ for all $\mu\in \Mon$, and (ii) $\mu \nu \preccurlyeq  \mu' \nu$ whenever $\mu, \mu', \nu \in \Mon$ are such that $\mu \preccurlyeq \mu'$. 
Let  $\prec$ be  any term order on $\Mon$ for which $x_0 \monord x_1 \monord \cdots \monord x_m$. For example, we can take $\prec$ to be the lexicographic order
$\lex$ defined by
\begin{equation}\label{def:lex}
x_0^{a_0} \cdots x_{m}^{a_{m} } \lex  x_0^{b_0} \cdots x_{m}^{b_{m} }  
\Longleftrightarrow 
\text{the first nonzero entry of }\mathrm{b}-\mathrm{a} \text{ is positive,}
\end{equation}
where $\mathrm{b}-\mathrm{a}$ denotes the difference tuple $(b_0 - a_0, \dots , b_m -a_m)$.
Other examples of such term orders are also possible. 
For now, we just fix a  term order $\prec$ on $\Mon$ for which $x_0 \monord x_1 \monord \cdots \monord x_m$. 
For a nonzero polynomial $F \in \Fq[x_0,\dots,x_m]$, we denote by $\lm_{\prec}(F)$ or simply $\lm(F)$, the \emph{leading monomial} of $F$, i.e., the largest monomial (w.r.t $\prec$) appearing in $F$ with a nonzero coefficient.
For any set $\s$ of nonzero polynomials in $\Fq[x_0,\dots,x_m]$ and any nonnegative integer $e$, we define
$$
\lm(\s):=\{\lm(F) : F \in \s\}, \quad \nabla_e(\s):=\nabla_e(\lm(\s)) \quad  \mathrm{and} \quad \Delta_e(\s):=\Delta_e(\lm(\s)).
$$
The following important result from \cite{BDG18} relates footprints in $\Monb_e$ to the number of $\Fq$-rational points of projective algebraic varieties defined over $\Fq$. 
Here, and hereafter, for an assertion depending on a nonnegative integer $e$, the expression
``for all $e \gg 0$" means that the assertion holds for all large enough values of $e$, i.e., there is a nonnegative integer $e_0$ such  that the assertion holds for all $e\ge e_0$.

\begin{theorem}[Projective $\Fq$-Footprint Bound] 
\label{thm:qfootprint}
Let $\s=\{F_1, \dots , F_r\}$ 
be a set of nonzero, projectively reduced homogeneous polynomials in $\Fq[x_0,\dots,x_m]$. Write $\x=V(\s)$ for the corresponding algebraic variety in $\PP^m$. Then 
$$
|\x(\Fq)| \le |\Delta_e(\s)| \quad \text{for all } e\gg 0.
$$
\end{theorem}

\begin{proof}
Since $F_1, \dots , F_r$ are projectively reduced, $\lm (\s) \subseteq \Monb$. Thus the
notion of footprint $\Delta_e(\s)$ defined above coincides with 
that of projective $\Fq$-footprint $\overline{\Delta}_e(\s)$ defined in \cite[Def. 3.10]{BDG18}. 
So the desired result follows from \cite[Thm. 3.12]{BDG18}.
\end{proof}
%
%

We can decompose 
$\Monb$ into disjoint subsets by considering the last variable appearing in a given monomial. Thus, following \cite{BDG18},
 we define
 $$
 \Monb^{(0)} = \{x_0^a : a \ge 0\} \quad \text{and  for $1 \le \ell \le m$,}\quad
 \Monb^{(\ell)}:=\{x_0^{a_0} \cdots x_{\ell}^{a_{\ell}} \in \Monb :  a_\ell >0\}.
 $$
Further, for any 
integers $e, \ell$ with $e\ge 0$ and $0\le \ell \le m$, define $\Monb_e^{(\ell)}:= \Monb^{(\ell)} \cap \Monb_{e}$.
Evidently, 
$\Monb = \coprod_{\ell=0}^m \Monb^{(\ell)}$ and $\Monb_e = \coprod_{\ell=0}^m \Monb_e^{(\ell)}$.
It is easy to see that 
$$
|\Monb_e^{(\ell)}|=q^\ell \ \text{ for } e \ge \ell(q-1)+1 \quad \text{and hence} \quad |\Monb_e|=\p_m \ \text{ for } e \ge m(q-1)+1.
$$
We refer to \cite[Section 2.1]{BDG18} for more details.

The sets $\Monb_e^{(\ell)}$ help us in decomposing the shadows or footprints into several disjoint components. To this end,
we define
\begin{equation}\label{eq:deltaeell}
\nabla_e^{(\ell)}(\s):=\nabla_e(\s) \cap \Monb_e^{(\ell)}  \quad \makebox{and} \quad \Delta_e^{(\ell)}(\s):=\Delta_e(\s) \cap \Monb_e^{(\ell)}.
\end{equation}
For any 
$\s \subseteq \Monb$ and any nonnegative integer $e$, it is clear that 
\begin{equation}\label{eq:disjointshfp}
|\nabla_e(\s)|= \sum_{\ell=0}^m |\nabla_e^{(\ell)}(\s)| \quad \makebox{and} \quad
|\Delta_e(\s)|= \sum_{\ell=0}^m |\Delta_e^{(\ell)}(\s)|.
\end{equation}

Given a nonnegative integer $\ell \le m$, by specializing the variables $x_{\ell+1},\dots,x_m$ to $1$, we
can associate to a subset of 
$\Monb$, a set of projectively reduced monomials in $x_0, \dots , x_{\ell}$ as follows.

\begin{definition}\label{def:sl}
Let $\s \subset \Monb$. For $0 \le \ell \le m$, we define
$$
\s^{\langle\ell\rangle}:=\{x_0^{a_0} \cdots x_{\ell}^{a_{\ell}}  \in \s : 0 \le a_j < q \ \mathrm{for \ all} \ j< \ell\}.$$
\end{definition}

Note  that if $0\le \ell \le m$ and if $\mu \in \s \setminus \s^{\langle\ell\rangle}$,  then either (i) $x_i \mid \mu$ for some $i > \ell$,  or (ii)  $\mu = x_0^{a_0} \cdots x_{j}^{a_{j}}$ with $a_j >q$ for some $j < \ell$.
In either case, it is easily seen that 
$\nabla_e^{(\ell)} (\mu) = \emptyset$. This shows that
\begin{equation}\label{eq:ssl}
\nabla_e^{(\ell)}(\s)=\nabla_e^{(\ell)}(\s^{\langle\ell\rangle}) \quad \text{and hence} \quad \Delta_e^{(\ell)}(\s)=\Delta_e^{(\ell)}(\s^{\langle\ell\rangle}) \quad \text{for any } e \ge 0.
\end{equation}

The following reformulation of Theorem \ref{thm:qfootprint} will be useful to us later.

\begin{corollary}\label{cor:qfootprint}
Let $F_1,\dots,F_r$ 
be any nonzero projectively reduced homogeneous polynomials in $\Fq[x_0,\dots,x_m]$,
and let  $\x=V(F_1,\dots,F_r)$ be the corresponding projective variety in $\PP^m$. 
Also, let $\s=\{\lm(F_1),\dots,\lm(F_r)\}.$ Then 
$$
|\x(\Fq)| \le \sum_{\ell=0}^m|\Delta_e^{(\ell)}(\s^{\langle\ell\rangle})| \quad \text{for all } e\gg 0.
$$
\end{corollary}

\begin{proof}
From equation \eqref{eq:ssl} and the second part of equation \eqref{eq:disjointshfp}, we  see 
that
\begin{equation}\label{eq:footprint}
|\Delta_e(\s)|= \sum_{\ell=0}^m |\Delta_e^{(\ell)}(\s^{\langle\ell\rangle})|.
\end{equation}
Thus the desired result follows from Theorem \ref{thm:qfootprint}.
\end{proof}

\section{Affine Combinatorics}
\label{combin}

In this section, we will consider some results from extremal combinatorics together with their generalizations and variants that will be useful for us later. These results are mainly concerned with sets of monomials in $\ell$ variables that are reduced in the usual (or the affine) sense. Throughout this section, $m$ will be a fixed positive integer and $\ell, d$ as well as $a_i, b_i$ denote nonnegative integers with $\ell \le m$.

\subsection{Shadows and Footprints in the Hypercube}
\label{subsec:3.1}
For $0 \le \ell \le m$, define
$$
\Nom^{(\ell)}:=\{x_0^{a_0}\cdots x_{\ell-1}^{a_{\ell-1}} \, : \, 0 \le a_j \le q-1  \makebox{ for all $j=0,\dots \ell-1$}\}.
$$
Note that the ``exponent map" given by
$x_0^{a_0}\cdots x_{\ell-1}^{a_{\ell-1}} \longmapsto (a_0,\dots,a_{\ell-1})$ lets us identify the set
$\Nom^{(\ell)}$ 
with  $Q^\ell$, where $Q:=\{0,1, \dots,q-1\}$, and thus we may refer to
$\Nom^{(\ell)}$ as the ($\ell$-dimensional) \emph{hypercube}. We remark that when $q=2$, elements of
$\Nom^{(\ell)}$ can be identified with subsets of $\{0, \dots , \ell-1\}$, whereas in general, they may be viewed as multisets formed by the elements  of $\{0, \dots , \ell-1\}$. We will, however, stick to viewing
$\Nom^{(\ell)}$ as the set of reduced monomials in $\ell$ variables $x_0, \dots , x_{\ell-1}$.

Divisibility of monomials gives a natural partial order on $\Nom^{(\ell)}$, which corresponds, via the exponent map, to the ``product order"  $\le_P$ on $Q^\ell$ defined
by  $(a_0,\dots,a_{\ell-1}) \le_P (b_0,\dots,b_{\ell-1})$ if and only if $a_i \le b_i$ for all $i=0,\dots,\ell-1$.
On the other hand, the usual lexicographic order on $Q^\ell$ corresponds to the total order on $\Nom^{(\ell)}$, which we denote, as in \eqref{def:lex}, by $\lex$. 
Note that 
if 
$\mu, \nu \in \Nom^{(\ell)}$ are such that 
%
$\nu \mid \mu$, i.e.,  $\nu$ divides $\mu$,  then $\nu \lexeq \mu$.

As before, for a nonnegative integer $d$, we define
$$
\Nom^{(\ell)}_d:=\{\mu \in \Nom^{(\ell)}\, : \, \deg \mu =d\} \quad \makebox{and} \quad
\Nom^{(\ell)}_{\le d}:=\{\mu \in \Nom^{(\ell)}\, : \, \deg \mu  \le d\}.
$$
The sets $\Nom^{(\ell)}_{< d}, \ \Nom^{(\ell)}_{\ge d}$ and $\Nom^{(\ell)}_{> d}$ are defined analogously. 

We will now define 
shadow and footprint in the context of the hypercube  $\Nom^{(\ell)}$. To avoid confusion with the notions defined in \S\,\ref{subsec:pfb}, we will use a different notation.

For any 
$\t \subseteq \Nom^{(\ell)}$, the \emph{shadow} and \emph{footprint} of $\t$ in $\Nom^{(\ell)}$ are denoted by $\SH^{(\ell)}(\t)$ and $\FP^{(\ell)}(\t)$, respectively, and defined by 
$$
\SH^{(\ell)}(\t):=\{\mu \in \Nom^{(\ell)} \, : \, \nu \mid \mu \ \makebox{for some $\nu \in \t$}\} \quad \makebox{and} \quad  \FP^{(\ell)}(\t):=\Nom^{(\ell)} \setminus \SH^{(\ell)}(\t).
$$
For a nonnegative integer $d$, we also define,
$$
\FP_d^{(\ell)}(\t):=\FP^{(\ell)}(\t) \cap \Nom^{(\ell)}_{d} \quad \makebox{and} \quad \FP_{\le d}^{(\ell)}(\t):=\FP^{(\ell)}(\t) \cap \Nom^{(\ell)}_{ \le d}.
$$
The sets $\FP^{(\ell)}_{< d}(\t), \ \FP^{(\ell)}_{\ge d}(\t)$ and $\FP^{(\ell)}_{> d}(\t)$ are defined 
analogously. 
Moreover, the corresponding subsets of $\SH^{(\ell)}(\t)$ are also defined in a similar manner.

\begin{definition}
Let $\ell, d, \rho, \rho'$ be integers satisfying $0\le d \le \ell(q-1)$, $0 \le \ell \le m$, $0 \le \rho \le |\Nom_{\le d}^{(\ell)}|$, and $0 \le \rho' \le |\Nom^{(\ell)}_d|.$
Define
\begin{eqnarray*}
\m_d^{(\ell)}(\rho) &: =& \text{the set of first $\rho$ elements of $\Nom_{\le d}^{(\ell)}$ in descending lexicographic order,} \\
\L_d^{(\ell)}(\rho') &: =& \text{the set of first $\rho'$ elements of $\Nom_{d}^{(\ell)}$ in descending lexicographic order}.
\end{eqnarray*}
\end{definition}

Note that if $\ell, \rho, \rho'$ are positive and $d<q$, then both $\m_d^{(\ell)}(\rho)$ and $\L_d^{(\ell)}(\rho')$ contain $x_0^d$, which is the largest element of $\Nom_{\le d}^{(\ell)}$ as well as $\Nom_{d}^{(\ell)}$ in 
lexicographic order. Also, 
for a fixed $d\ge 0$, the set $\Nom_{\le d}^{(\ell)}$ is finite and $\lex$ is a total order on it.
Consequently,  a 
set of the form $\m_d^{(\ell)}(\rho)$ can be characterized as a
subset $\t$ of $\Nom_{\le d}^{(\ell)}$ that is \emph{upwards closed}, which means $\mu \in \t$ whenever
$\mu, \mu' \in \Nom_{\le d}^{(\ell)}$ with $\mu' \lex \mu$ and $\mu' \in \t$. 
Similarly, sets of the form $\L_d^{(\ell)}(\rho')$ can be characterized as upwards closed subsets of
$\Nom_{d}^{(\ell)}$.

\subsection{Extremal Combinatorics}
We will now discuss some combinatorial results that help us determine the subsets of $\Nom_{d}^{(\ell)}$ (resp. $\Nom_{\le d}^{(\ell)}$) of a given cardinality that have footprint of the maximum possible size.   
We begin with a result due to Heijnen and Pellikaan \cite[Prop. 5.9]{HP}. Its formulation below is as in \cite[Lem. 4.2]{BD}, where  the result
is proved in a more general setting.

\begin{proposition}
\label{prop:shadow}
Let $\ell, d, \rho$ be integers satisfying $1\le \ell \le m$, $0\le d \le \ell(q-1)$, and $1\le \rho \le |\Nom_{\le d}^{(\ell)}|$. Also,
let $\alpha$ be the $\rho^{\rm th}$ element of   $\Nom_{\le d}^{(\ell)}$, i.e., the smallest element of
$\m_d^{(\ell)}(\rho)$ in lexicographic order. Then
$$
\SH^{(\ell)}(\m_d^{(\ell)}(\rho))=\{\mu \in \Nom^{(\ell)} \, : \, \alpha \lexeq \mu\}.
$$
Consequently,  $\SH^{(\ell)}(\m_d^{(\ell)}(\rho)) \cap \Nom_{\le d}^{(\ell)} = \m_d^{(\ell)}(\rho)$.
\end{proposition}


\begin{corollary}\label{cor:compressed}
Let $\ell, d, \rho$ be integers satisfying $1\le \ell \le m$, $1\le d \le \ell(q-1)$, and $0 \le \rho \le |\Nom_{\le d}^{(\ell)}|$. Then $\SH_d^{(\ell)}(\m_{d-1}^{(\ell)}(\rho))=\L_d^{(\ell)}(\rho')$ for some
nonnegative  integer $\rho'$. 
Moreover, if $\rho$ is positive, then so is $\rho'$.
\end{corollary}

\begin{proof}
If $\rho = 0$, then $\m_{d-1}^{(\ell)}(\rho)$ is empty, and hence so is its shadow in $\Nom^{(\ell)}$.
Thus we can take $\rho'=0$ in this case.
Now suppose $\rho \ge 1$. 
First, note that
$\SH_d^{(\ell)}(\m_{d-1}^{(\ell)}(\rho))$ is nonempty. 
Indeed,  $d-1<\ell(q-1)$ and so we can write $d-1 = j(q-1)+a_j$ for unique integers $j, a_j$ with
$0\le j  \le \ell -1$ and $0\le a_j < q-1$. Since $\rho\ge 1$, the set   $\m_{d-1}^{(\ell)}(\rho)$ contains $\nu:= x_0^{q-1}\cdots x_{j-1}^{q-1}x_j^{a_j}$, 
being the largest element of  this set 
in lexicographic order. Now $\mu:= \nu x_{\ell -1}$ is clearly in $\SH_d^{(\ell)}(\m_{d-1}^{(\ell)}(\rho))$. 

To complete the proof, it suffices to show that $\SH_d^{(\ell)}(\m_{d-1}^{(\ell)}(\rho))$ is upwards closed. Assume the contrary.
Then there exist $\mu, \mu' \in \Nom_d^{(\ell)}$ such that $\mu' \in \SH_d^{(\ell)}(\m_{d-1}^{(\ell)}(\rho))$ and
$\mu' \lex \mu$,  but $\mu \not\in \SH_d^{(\ell)}(\m_{d-1}^{(\ell)}(\rho))$. 
Let $\nu \in \m_{d-1}^{(\ell)}(\rho)$ be such that $\nu \mid \mu'.$  Then $\nu \lexeq \mu'$.
Also, let $\nu'$ be the least element of $\m_{d-1}^{(\ell)}(\rho)$ in lexicographic order.
Then $\nu' \lexeq \nu$ and $\nu \lexeq \mu'$, since $\nu \mid \mu'$. Consequently, 
$\nu' \lex \mu.$ By the previous proposition, $\mu \in \SH_d^{(\ell)}(\m_{d-1}^{(\ell)}(\rho)),$ which is a contradiction.
\end{proof}

The following 
result is due to Clements and Lindstr\"om \cite[Cor. 1]{CL} (see also \cite[Thm. 3.1]{BD}). 
The case $q=2$ of it 
is equivalent to the Kruskal-Katona theorem.

\begin{theorem}[Clements--Lindstr\"om]
\label{thm:CL}
Let $\t \subseteq \Nom_d^{(\ell)}$ with $|\t| = \rho'$. 
Then
$$
\SH_{d+1}^{(\ell)}(\L_d^{(\ell)}(\rho')) \subseteq \L_{d+1}^{(\ell)}(|\SH_{d+1}^{(\ell)}(\t)|) \ \mathrm{and \ hence} \ |\FP_{d+1}^{(\ell)}(\t)| \le |\FP_{d+1}^{(\ell)}(\L_d^{(\ell)}(\rho'))| .
$$
\end{theorem}

As in \cite[Cor. 3.2]{BD}, the following corollary can be deduced easily from the above theorem 
by using induction on $e$.

\begin{corollary}\label{cor:Wei}
Let $\t \subseteq \Nom_d^{(\ell)}$ with $|\t| = \rho'$, and let $e$ be an integer  $\ge d$.
Then
$$
\SH_{e}^{(\ell)}(\L_d^{(\ell)}(\rho')) \subseteq \L_{e}^{(\ell)}(|\SH_{e}^{(\ell)}(\t)|); \ \text{ in particular, } \  \left| \SH_{e}^{(\ell)}(\L_d^{(\ell)}(\rho')) \right| \le |\SH_{e}^{(\ell)}(\t)|.
$$
Consequently,
$|\FP_{e}^{(\ell)}(\t)| \le |\FP_{e}^{(\ell)}(\L_d^{(\ell)}(\rho'))|$ and hence 
$|\FP^{(\ell)}(\t)| \le |\FP^{(\ell)}(\L_d^{(\ell)}(\rho'))|.$
\end{corollary}

The following theorem can be traced back to Wei \cite[Lemma 6]{W} and it
gives an analogue of the last inequality for subsets of $\Nom_{\le d}^{(\ell)}$.

\begin{theorem}[Wei]  \label{thm:Wei}
Let $\t \subset \Nom_{\le d}^{(\ell)}$ be a subset with $|\t|=\rho.$ Then 
$$ |\FP^{(\ell)}(\t)| \le |\FP^{(\ell)}(\m_d^{(\ell)}(\rho))|.$$
\end{theorem}

Strictly speaking, Lemma 6 of Wei \cite{W} proves the above theorem for the special case when $q=2$. 
A general version is stated in Heijnen and Pellikaan \cite[Thm.~5.7]{HP}, although in a slightly different way. A detailed proof appears in Appendix A of Heijnen's thesis \cite{H}. We refer to \cite[Thm. 3.8]{BD} for a similar result in a more general setting. We conclude this subsection by proving the  following common generalization of the results of Clements--Lindstr\"om (Corollary \ref{cor:Wei}) and of Wei and Heijnen--Pellikaan  (\hbox{Theorem \ref{thm:Wei}}).  Indeed, Corollary \ref{cor:Wei} corresponds to the case $\rho'=\rho$, while Theorem \ref{thm:Wei} corresponds to the case $\rho'= 0$.

\begin{theorem}\label{thm:affinecomb}
Assume that $\ell, d, \rho$ are positive integers with $\ell \le m$,  $d\le \ell (q-1)$, and  $\rho \le |\Nom_{\le d}^{(\ell)}|$.
Let $\t \subseteq \Nom^{(\ell)}_{\le d}$ with $|\t| = \rho$. If $ \rho' := |\t \cap \Nom_d^{(\ell)}|$, then
$$
|\FP^{(\ell)}(\t)| \le |\FP^{(\ell)}(\u)|, \quad \makebox{where} \quad \u:=\L_d^{(\ell)}(\rho') \cup \m_{d-1}^{(\ell)}(\rho - \rho').
$$
\end{theorem}

\begin{proof}
Let $ \rho' := |\t \cap \Nom_d^{(\ell)}|$ and $ \u:=\L_d^{(\ell)}(\rho') \cup \m_{d-1}^{(\ell)}(\rho - \rho')$.
It suffices to show that $|\SH^{(\ell)}(\u)| \le |\SH^{(\ell)}(\t)|.$
We will do this by distinguishing two cases.

\smallskip

\noindent
{\bf Case 1:} $\L^{(\ell)}_d (\rho') \subseteq \SH^{(\ell)}_d (\m_{d-1}^{(\ell)}(\rho-\rho'))$.

Let $\mu \in \SH^{(\ell)}(\u)$. Then there is $\nu \in \u$ such that $\nu \mid \mu$.  Suppose
$\nu \in \L^{(\ell)}_d (\rho') $.  Since $\L^{(\ell)}_d (\rho') \subseteq \SH^{(\ell)}_d (\m_{d-1}^{(\ell)}(\rho-\rho'))$, there is $\nu' \in \m_{d-1}^{(\ell)}(\rho-\rho')$ such that
$\nu' \mid \nu$, and hence $\nu'\mid \mu$. 
This shows that  $\SH^{(\ell)}(\u) = \SH^{(\ell)} ( \m_{d-1}^{(\ell)}(\rho-\rho') )$.
Consequently,
$$
 |\SH^{(\ell)}(\u)\big| = \big|\SH^{(\ell)}(\m_{d-1}^{(\ell)}(\rho-\rho'))\big| \le \big|\SH^{(\ell)}(\t \cap \Nom_{\le d-1}^{(\ell)})\big| \le \big| \SH^{(\ell)}(\t) \big|,
$$
where the penultimate inequality is a consequence of Theorem \ref{thm:Wei} (applied to shadows instead of footprints), while the last inequality follows since $\t \cap \Nom_{\le d-1}^{(\ell)} \subseteq \t$.

\smallskip

\noindent
{\bf Case 2:} $\L_d^{(\ell)}(\rho') \not\subseteq \SH^{(\ell)}_d(\m_{d-1}^{(\ell)}(\rho-\rho'))$.

\smallskip

By Corollary \ref{cor:compressed},  
$\SH^{(\ell)}_d(\m_{d-1}^{(\ell)}(\rho-\rho')) = \L_d^{(\ell)}(\rho'')$ for some nonnegative integer $\rho''$. Hence $\L_d^{(\ell)}(\rho') \not\subseteq \L_d^{(\ell)}(\rho'')$ and this implies that $\rho'\not\le \rho''$, i.e., $\rho'' < \rho'$. It follows that $\SH^{(\ell)}_d(\m_{d-1}^{(\ell)}(\rho-\rho')) \subseteq \L_d^{(\ell)}(\rho')$. Consequently,
\begin{equation}\label{eq:shadow2}
\SH_{\ge d}^{(\ell)}(\u)
= \SH_{\ge d}^{(\ell)}(\L_d^{(\ell)}(\rho')) \cup \SH_{\ge d}^{(\ell)}(\m_{d-1}^{(\ell)}(\rho-\rho'))
= \SH_{\ge d}^{(\ell)}(\L_d^{(\ell)}(\rho')).
\end{equation}
On the other hand, $\SH^{(\ell)}_{<d} (\u)=\SH^{(\ell)}_{<d} (\m_{d-1}^{(\ell)}(\rho-\rho'))=\SH^{(\ell)}(\m_{d-1}^{(\ell)}(\rho-\rho')) \cap \Nom_{< d}^{(\ell)}$, and so by Proposition \ref{prop:shadow}.
we see that
\begin{equation}\label{eq:shadow1}
\SH^{(\ell)}_{<d} (\u)= 
\m_{d-1}^{(\ell)}(\rho-\rho').
\end{equation}
Using equations \eqref{eq:shadow2} and \eqref{eq:shadow1}, we obtain
$$
|\SH^{(\ell)} (\u)|  
= \rho - \rho'+|\SH^{(\ell)}_{\ge d} (\L_d^{(\ell)}(\rho'))|  \le \rho-\rho'+|\SH^{(\ell)}_{\ge d} (\t \cap \Nom^{(\ell)}_d)|,
$$
where the last inequality follows from Corollary \ref{cor:Wei}. Now since $|\t \cap \Nom^{(\ell)}_{<d}| = \rho - \rho'$, 
$$
\big|\SH^{(\ell)} (\u) \big|  \le \big|\SH^{(\ell)}_{<d} (\t \cap \Nom^{(\ell)}_{<d})\big| + \big|\SH^{(\ell)}_{\ge d} (\t \cap \Nom^{(\ell)}_d)\big|.
$$
Hence $|\SH^{(\ell)} (\u)| \le |\SH^{(\ell)}_{<d} (\t)| +|\SH^{(\ell)}_{\ge d} (\t)|= |\SH^{(\ell)} (\t)|$, as desired.
\end{proof}

\section{Specializations and Expanders}
\label{combin2}
%
In order to effectively relate the two notions $\Delta_e$ and $\FP$ of footprint considered in the previous two sections, we will introduce two maps, denoted $\f^{(\ell)}$ and $\phi$, on the space of projectively reduced monomials in $x_0, \dots , x_{m}$ and prove some of their properties. Throughout this section, $m$ is a fixed positive integer, while $\ell, d, e$ denote nonnegative integers satisfying
$\ell \le m$.

\subsection{Specialization}
For any nonnegative integer $\ell \le m$, we define 
$$
\f^{(\ell)}: \Monb \to \Monb \cup \{0\}
\quad \makebox{by} \quad \f^{(\ell)}(\mu) := 
\begin{cases} 
x_0^{a_0}\cdots x_{\ell-1}^{a_{\ell-1}} & \makebox{if }\mu =x_0^{a_0}\cdots x_{\ell-1}^{a_{\ell-1}}x_\ell^{a_\ell},\\
0 & \makebox{if } \mu \not\in \Fq[x_0,\dots,x_\ell], \\
\end{cases} 
$$
with the usual convention that an empty product equals $1$. 
We may refer to $\f^{(\ell)}$ as the \emph{specialization map} at level $\ell$, since it corresponds to specializing the variables $(x_\ell,x_{\ell+1},\dots,x_m)$ to $(1,0,\dots,0)$.

\begin{proposition}\label{prop:divisible}
Let $\ell, d, e$ be nonnegative integers such that $e \ge d+m(q-1)$ and $\ell \le m$.
Suppose $\mu  \in \Monb_d\cap \Fq[x_0,\dots,x_\ell]$ and $\nu \in \Monb_e^{(\ell)}.$
Then 
$$
\mu \mid \nu  \iff \f^{(\ell)}(\mu) \mid \f^{(\ell)}(\nu).
$$
\end{proposition}

\begin{proof}
We can write $\mu = x_0^{a_0} \cdots x_{\ell}^{a_{\ell}}$ 
and $\nu=x_0^{b_0}\cdots x_{\ell-1}^{b_{\ell-1}}x_\ell^{b_\ell}$ for some nonnegative integers $a_j, b_j$ for $0\le j \le \ell$ such that
$b_j \le q-1$
for all $j \le \ell-1$ and $b_{\ell} >0$.  Note that
$\f^{(\ell)}(\mu)=x_0^{a_0}\cdots x_{\ell-1}^{a_{\ell-1}}$ and $\f^{(\ell)}(\nu)=x_0^{b_0}\cdots x_{\ell-1}^{b_{\ell-1}}.$

If $\mu \mid \nu$, then $a_i \le b_i$ for $0\le i \le \ell$, and this 
readily implies that $\f^{(\ell)}(\mu) \mid \f^{(\ell)}(\nu).$ To prove the converse, suppose  $\f^{(\ell)}(\mu) \mid \f^{(\ell)}(\nu).$ Then $a_i \le b_i \le q-1$ for 
$0 \le i \le \ell-1.$
Since $\sum_{i=0}^{\ell } a_i = d$ and $\sum_{i=0}^{\ell} b_i = e$, we obtain
$$
b_{\ell} - a_{\ell} = e-d + \sum_{i=0}^{\ell - 1} (a_i - b_i) \ge m(q-1) - \ell  (q-1)  \ge 0.
$$
This shows that $\mu \mid \nu$. 
\end{proof}

\begin{remark}
\label{rem:sigmam}
Suppose $0\le d < q$. Then $\sigma^{(m)}$ gives a bijection of $\Monb_d = \Mon_d$ onto $\Nom^{(m)}_{\le d}$. Indeed, if $\mu = x_0^{a_0} \cdots x_{m}^{a_{m}} \in \Monb_d$, then clearly,  $0\le a_i \le d \le q-1$  for $0\le i \le m$; 
also $a_0 + \dots + a_{m -1} \le d$. Moreover, $a_m = d -  a_0 - \dots - a_{m -1}$, and so $\mu$ is determined by $\f^{(m)}(\mu)=x_0^{a_0}\cdots x_{m-1}^{a_{m-1}} \in \Nom^{(m)}_{\le d}$.
A similar reasoning shows that $\sigma^{(m)}$ preserves lexicographic order, i.e., for any
$\mu, \nu \in \Monb_d$, we have $\mu \lex \nu \iff \sigma^{(m)}(\mu) \lex \sigma^{(m)}(\nu)$.
\end{remark}

The following result gives a useful relation between the two notions of footprint.

\begin{theorem}\label{thm:projectivetoaffine}
Let $d, e$ be any nonnegative integers 
such that  $e \ge d + m(q-1)$ and let $\s \subset \Monb_d$. Then
$$
|\Delta_e^{(\ell)}(\s^{\langle\ell\rangle})|=|\FP^{(\ell)}(\f^{(\ell)}(\s^{\langle\ell\rangle}))| \quad \text{for all nonnegative integers } \ell \le m, 
$$
and consequently,
$$
|\Delta_e(\s)|=\sum_{\ell=0}^{m}|\FP^{(\ell)}(\f^{(\ell)}(\s^{\langle\ell\rangle}))|.
$$
\end{theorem}

\begin{proof}
Fix a nonnegative integer $\ell \le m$. It is clear from Definition \ref{def:sl}  that $\f^{(\ell)}(\s^{\langle\ell\rangle}) \subseteq \Nom^{(\ell)}.$
From Proposition \ref{prop:divisible}, we see that the map
$$
\Delta_e^{(\ell)}(\s^{\langle\ell\rangle}) \to \FP^{(\ell)}(\f^{(\ell)}(\s^{\langle\ell\rangle})) \quad \text{defined by} \quad \mu \longmapsto \f^{(\ell)}(\mu)
$$
is well-defined. Moreover, since $e \ge d + m(q-1)$, this map is easily seen to be a bijection. This yields the first assertion in the theorem. Consequently, we obtain the last assertion from equation  \eqref{eq:footprint}.
\end{proof}

\subsection{Footprint Expander}
In this subsection, we consider a degree-preserving map $\phi$ on sets of projectively reduced monomials in $x_0, \dots , x_m$ such that 
$\phi$ is injective and has the property that
$|\Delta_e(\s)|\le |\Delta_e(\phi(\s))|$ for any $\s \subseteq \Monb$ and $e\gg 0$.  For this reason, $\phi$ may be referred to as an \emph{expander} map.

\begin{definition}\label{def:phi}
Let $\s \subseteq \Monb$. 
If $\mu=x_0^{i_0}\cdots x_m^{i_m} \in \s$, then define
$$
\phi(\mu) :=
\begin{cases}
\mu  \frac{x_{m-1}}{x_m}  &  \mathrm{if} \   x_0^{i_0}\cdots x_{m-2}^{i_{m-2}}x_{m-1}^{i_{m-1}+i_m} \not\in \s  \ \mathrm{and} \  i_{m-1} + 1 <q,\\
\mu  &  \mathrm{otherwise.}
\end{cases}
$$
\end{definition}
Note that if $i_m =0$, then  $x_0^{i_0}\cdots x_{m-2}^{i_{m-2}}x_{m-1}^{i_{m-1}+i_m} = \mu \in \s$. Hence, in the first case of the definition, we must have 
$i_m >0$. In particular, $\phi (\mu)$ is always a monomial, and we obtain a well-defined map $\phi: \s \to \Monb$, which preserves degrees.


\begin{proposition}
Let $\s \subseteq \Monb$. Then the map $\phi: \s \to \Monb$ is injective and it satisfies
$$
\s^{\langle m-1 \rangle}=\phi(\s^{\langle m-1 \rangle}) \subseteq \phi(\s)^{\langle m-1 \rangle}
\quad \text{and} \quad \phi(\s^{\langle m \rangle}) = \phi(\s)^{\langle m \rangle}.
$$
\end{proposition}

\begin{proof}
It is easy to see that $\phi$ is injective. Also, as noted earlier, $\phi(\mu) = \mu$ in case $\mu \in \s^{\langle m-1 \rangle}$. This implies that $\s^{\langle m-1 \rangle}=\phi(\s^{\langle m-1 \rangle}) \subseteq \phi(\s)^{\langle m-1 \rangle}$. To prove that $\phi(\s^{\langle m \rangle}) = \phi(\s)^{\langle m \rangle}$, suppose   $\mu = x_0^{i_0} \cdots x_m^{i_m} \in \s^{\langle m \rangle}$. 
Then $i_j < q$ for $0 \le j \le m-1$.
In case $\phi(\mu) = \mu$, then clearly,  $\phi(\mu) \in \phi(\s)^{\langle m \rangle}$. In particular, if $i_{m-1} = q-1$, then $\phi (\mu) \in \phi(\s)^{\langle m \rangle}$ because in this case, $i_{m-1} + 1 = q$ and so $\phi(\mu) = \mu$.  On the other hand, if $i_{m-1} < q-1$ and $\phi (\mu) \neq \mu$, then $\phi(\mu) = x_0^{i_0} \cdots x_{m-1}^{i_{m-1} + 1} x_m^{i_m -1}$. Since $i_{m-1} + 1 < q$, we obtain $\phi(\mu) \in \phi(\s)^{\langle m \rangle}$. This shows that $\phi(\s^{\langle m \rangle}) \subseteq \phi(\s)^{\langle m \rangle}$. In order to prove the reverse inclusion, we take $\mu \in \phi(\s)^{\langle m \rangle}$. Since $\mu \in \phi (\s)$, there exists $\mu' \in \s$ such that $\phi(\mu') = \mu$. It is trivial to see that $\mu' \in \s^{\langle m \rangle}$. 
\end{proof}

We now give a series of lemmas and a proposition leading to the result that $\phi$ does not decrease footprints in $\Monb_e$ for all large enough $e$.

\begin{lemma}\label{claim1}
Let $\s$ be a finite subset of $\Monb$. Then $\Delta_e^{(m)}(\s) \subseteq \Delta_e^{(m)}(\phi(\s))$ for all $e\gg 0$.
\end{lemma}

\begin{proof}
It suffices to show that $\nabla_e^{(m)} (\phi(\s)) \subseteq \nabla_e^{(m)} (\s)$ for all $e\gg 0$.
To this end, let $\nu \in \nabla_e^{(m)} (\phi(\s)).$ Then $\nu \in \Monb_e^{(m)}$ and there exists $\mu \in \phi(\s)$ such that $\mu \mid \nu.$
If $\mu \in \s$, then $\nu \in \nabla_e^{(m)} (\s).$
Suppose $\mu \not \in \s$.  Since $\mu \in \phi(\s)$, it follows~from Definition~\ref{def:phi} that
$\mu x_m/x_{m-1} \in \s.$ In particular, $x_{m-1} \mid \mu$ and $(\mu/x_{m-1}) \mid \nu.$
Further, since $\nu \in \Monb_e^{(m)}$, we see that $\deg_{x_m} \nu \gg 0$ for all $e\gg 0$, whereas
$\deg_{x_m} \mu$ is bounded since $\s$ is finite. Consequently, 
$(\mu x_m/x_{m-1} )\mid \nu$,  
and so  $\nu \in \nabla_e^{(m)} (\s)$ for 
$e\gg 0$.
\end{proof}

\begin{lemma}\label{claim2}
Let $\s \subseteq \Monb$. 
Then $\Delta_e^{(m-1)}(\phi(\s)) \subseteq \Delta_e^{(m-1)}(\s)$ for all $e\ge 0$.
\end{lemma}
\begin{proof}
Fix a nonnegative integer $e$.  
We will show that $\nabla_e^{(m-1)} (\s) \subseteq \nabla_e^{(m-1)} (\phi(\s))$.
Let $\nu \in \nabla_e^{(m-1)} (\s)$. Then $\nu \in \Monb_e^{(m-1)}$ and there exists $\mu \in \s$ such that $\mu \mid \nu.$ In particular, $\deg_{x_m} \mu \le \deg_{x_m} \nu = 0$, and so
$\phi(\mu) = \mu$. It follows that $\mu \in \phi(\s)$ and so $\nu \in \nabla_e^{(m-1)} (\phi(\s)).$
%
\end{proof}

\begin{lemma}\label{claim3}
Let $\s \subseteq \Monb$ and  let  $e$ be a nonnegative integer.  Suppose $\nu \in \Delta_e^{(m-1)}(\s)$ and $\mu \in \Monb$ satisfy $\mu \mid \nu$. Then $\mu \not \in \s^{\langle m-1 \rangle}.$
\end{lemma}

\begin{proof}
Since $\Delta_e^{(m-1)}(\s) = \Delta_e^{(m-1)}(\s^{\langle m-1 \rangle})$, this follows directly.
\end{proof}

\begin{lemma}\label{nonempty}
Let $\s \subseteq \Monb$ and let $e$ be a nonnegative integer. Suppose there exists
$\nu \in \Delta_e^{(m-1)}(\s) \setminus \Delta_e^{(m-1)}(\phi(\s))$, i.e., 
$\nu \in \nabla_e^{(m-1)}(\phi(\s)) \setminus \nabla_e^{(m-1)}(\s).$ Then the set
$\{ \mu \in \s :  \deg_{x_{m-1}} \mu <q \ \makebox{and} \ \mu \mid (\nu  x_m)\}$
is nonempty.
\end{lemma}
\begin{proof}
Since $\nu \in \nabla_e^{(m-1)}(\phi(\s)) = \nabla_e^{(m-1)}(\phi(\s)^{\langle m-1 \rangle})$, there exists $\tilde{\mu} \in \phi(\s)^{\langle m-1 \rangle}$ such that $\tilde{\mu} \mid \nu$. Moreover, $\nu \not\in \nabla_e^{(m-1)}(\s)$ implies $\tilde{\mu} \not\in \s.$ Hence $x_{m-1} \mid \tilde{\mu}$ and $\tilde{\mu}x_m/x_{m-1} \in \s.$ 
Let $\mu := \tilde{\mu}x_m/x_{m-1}$.
Note that $\mu$ is projectively reduced since $\mu \in \s$.
Further, $\deg_{x_m} \mu > 0$ since $\phi(\mu) \ne \mu$. This implies that $\deg_{x_{m-1}} \mu <q.$
Also, since $\tilde{\mu} \mid \nu$, we see that $\tilde{\mu}/x_{m-1}$ divides $\nu$ and so 
$\mu \mid (\nu x_m)$.
\end{proof}

\begin{proposition}\label{injective}
Let $\s \subset \Monb$ be a finite set. Then 
$$
\big|\Delta_e^{(m-1)}(\s) \setminus \Delta_e^{(m-1)}(\phi(\s))\big| \le \big|\Delta_e^{(m)}(\phi(\s)) \setminus \Delta_e^{(m)}(\s) \big| \quad \text{for all $e\gg 0$. }
$$
\end{proposition}
\begin{proof}
For a nonnegative integer $e$ and $\nu  \in \Delta_e^{(m-1)}(\s) \setminus \Delta_e^{(m-1)}(\phi(\s))$,  
define 
$$
S_\nu:=\{ \mu \in \s :  \deg_{x_{m-1}} \mu <q \text{ and } \mu \mid (\nu  x_m^j) \ \mathrm{ for \ some} \ j >0  \}.
$$
By Lemma \ref{nonempty},  
the set  $S_\nu$ is nonempty. Moreover, since $\s$ is a finite set, we see that 
$$
k_{\nu} := \min \{ \deg_{x_{m-1}}\! \mu : \mu \in S_{\nu}\} \quad \text{and} \quad
E:= \max \{\deg_{x_m} \! \mu : \mu \in \s\}
$$
are well-defined, and
$$
S_\nu =\{ \mu \in \s :  \deg_{x_{m-1}} \mu <q \text{ and } \mu \mid (\nu  x_m^E ) \}.
$$
For $e \ge 0$, consider the map
$
\psi: \Delta_e^{(m-1)}(\s) \setminus \Delta_e^{(m-1)}(\phi(\s)) \to \Monb_e
$ 
given by
$$
\nu=x_0^{i_0}\cdots x_{m-1}^{i_{m-1}} \longmapsto x_0^{i_0}\cdots x_{m-2}^{i_{m-2}} x_{m-1}^{k_\nu}x_{m}^{i_{m-1}-k_\nu}.
$$
Clearly, the map $\psi$ is injective and to prove the proposition, it is enough to show that the image of $\psi$ is contained in $\Delta_e^{(m)}(\phi(\s)) \setminus \Delta_e^{(m)}(\s)$, provided $e\gg 0$. 

Let us fix $e\ge 0$ and $\nu = x_0^{i_0}\cdots x_{m-1}^{i_{m-1}} \in \Delta_e^{(m-1)}(\s) \setminus \Delta_e^{(m-1)}(\phi(\s))$. Also, for simplicity, let us write $k= k_{\nu}$. From Lemma \ref{nonempty}, we see that 
$k <q$.
By the definition of $k$, there exists $\mu^* \in \s$ with $\deg_{x_{m-1}} \mu^*=k$ such that $\mu^* \mid (\nu x_m^E)$.
Now since $\nu \in \Monb^{(m-1)}_e$, we see that $i_j < q$ for $0\le j < m-1$, whereas $i_{m-1} \gg 0$ if $e\gg 0$. Hence $E \le i_{m-1} - k$ if $e \gg 0$.  This implies that
$\mu^* \mid \psi(\nu)$ whenever  $e \gg 0$. Since $\mu^* \in \s$, we conclude that $\psi(\nu) \not \in \Delta_e^{(m)}(\s).$
We now prove that $\psi(\nu) \in \Delta_e^{(m)}(\phi(\s))$ for $e\gg 0$ by distinguishing  two cases.

\smallskip

\noindent
\textbf{Case 1:} $k=q-1.$

Since $\nu \in \nabla_e^{(m-1)}(\phi(\s))$, there exists $\mu \in \phi(\s)$ such that $\mu \mid \nu$.
Let $\tilde{\mu} \in \s$ be such that $\mu=\phi(\tilde{\mu}).$ 
Then either $\tilde{\mu}=\mu$ or $\tilde{\mu}=\mu x_m/x_{m-1}.$ If $\tilde{\mu} = \mu$, then $\nu \in \nabla_e^{(m-1)}(\s)$, which is a contradiction. Thus $\tilde{\mu}=\mu x_m/x_{m-1}.$
Since $\mu \mid \nu,$ clearly $\tilde{\mu} \mid (\nu x_m)$. Hence $\tilde{\mu} \in S_\nu.$ From the assumption that $k=q-1$, we conclude that $\deg_{x_{m-1}}(\tilde{\mu})=q-1.$
But then $\mu=\phi(\tilde{\mu})=\tilde{\mu}$ which is again 
a contradiction. Therefore, the case $k = q-1$ can not occur.

\smallskip

\noindent
\textbf{Case 2:} $k<q-1.$

Suppose, if possible,  $\psi(\nu) \in  \nabla_e^{(m)}(\phi(\s))$. Then there  exists
$\mu \in \phi(\s)$ such that $\mu \mid \psi(\nu)$. Write $\mu=x_0^{j_0}\cdots x_{m}^{j_{m}}$.  Then $j_0 \le i_0, \dots, j_{m-2} \le i_{m-2}$ and $j_{m-1} \le k <q-1$.
Note that  $x_0^{j_0}\cdots x_{m-2}^{j_{m-2}} x_{m-1}^{j_{m-1}+j_m} \mid \nu$ whenever $e \gg 0$.  This implies that $x_0^{j_0}\cdots x_{m-2}^{j_{m-2}} x_{m-1}^{j_{m-1}+j_m} \not \in \s,$ since $\nu \in \Delta_e^{(m-1)}(\s).$
Choose $\tilde{\mu} \in \s$ such that $\phi (\tilde{\mu}) = \mu$. If $\tilde{\mu} = \mu$, then  by the  definition of $\phi$ and the fact that  $x_0^{j_0}\cdots x_{m-2}^{j_{m-2}} x_{m-1}^{j_{m-1}+j_m} \not \in \s$, we conclude that $j_{m-1} \ge q-1$, which is a contradiction.
Therefore 
$\tilde{\mu} \neq \mu.$ This implies that  $\tilde{\mu} = x_0^{j_0}\cdots x_{m-1}^{j_{m-1}-1}x_{m}^{j_{m}+1} \in \s.$ Furthermore, $\tilde{\mu} \mid \nu x_m^j$ whenever $j \ge j_m + 1$ and $e \gg 0$.
Thus $\tilde{\mu} \in S_{\nu}$. 
But then $\deg_{x_{m-1}} \tilde{\mu} = j_{m-1} - 1 < j_{m-1} \le k$, and this contradicts the minimality of $k$ if $e\gg 0$.  Hence $\psi(\nu) \in  \Delta_e^{(m)}(\phi(\s))$ 
for  $e\gg 0$. 
\end{proof}

\begin{theorem}\label{phi} Let $\s \subset \Monb$ be a finite set. Then
$|\Delta_e(\s)| \le |\Delta_e(\phi(\s))|$ for 
$e \gg 0$.
\end{theorem}

\begin{proof}
Since $\phi(\mu)=\mu$ for all $\mu \in \s$ with $\deg_{x_m} \mu=0$ and since $\deg_{x_{m-1}} \phi (\mu) > 0$ if $\phi (\mu) \neq \mu$, it is clear that
$\s^{\langle j \rangle} = \phi(\s)^{\langle j \rangle}$ for $j=0, \dots, m-2$. Hence
$$\Delta_e^{(j)} (\s) = \Delta_e^{(j)} (\s^{\langle j \rangle}) = \Delta_e^{(j)} (\phi(\s)^{\langle j \rangle})= \Delta_e^{(j)} (\phi(\s)) \ \ \mathrm{for \ all} \ j =0, \dots, m-2.$$
Thus it is enough to show that
$$|\Delta_e^{(m)}(\s)| + |\Delta_e^{(m-1)}(\s)| \le |\Delta_e^{(m)}(\phi (\s))| + |\Delta_e^{(m-1)}(\phi (\s))| .$$ This follows directly from Proposition \ref{injective} in view of
Lemmas \ref{claim1} and \ref{claim2}. 
\end{proof}

\section{Number of Solutions of Equations over Finite Fields}
\label{sec:er}

Throughout this section, $m$ denotes a fixed positive integer, while $d, r$ are nonnegative integers. 
Recall that $e_r(d,m)$ has been defined by \eqref{erdm} and the conventions \eqref{convention-erdm} whenever 
$ r \le \binom{m+d}{d}$. 
In this section, we shall prove our main results concerning 
$e_r(d,m)$
and a related quantity, called $\overline{e}_r(d,m)$ that we shall define shortly.

\subsection{Projectively Reduced Equations}\label{subsec:RedEqns}
As in \S\,\ref{subsec:projred},
let $\Gamma_q(\Fq)$ denote  the ideal of $\Fq[x_0, \dots, x_m]$ generated by
$\{x_i^qx_j-x_ix_j^q \mid 0 \le i < j \le m\}$. 
Evidently, $\Gamma_q(\Fq)$ is a homogeneous ideal and if we let $\Gamma_q(\Fq)_d$ denote its $d^{\rm th}$ homogeneous component, then its vector space dimension is known, e.g.,  from
\cite[Thm. 5.2]{MR}, namely 
\begin{equation}\label{eq:rd}
r_d:= \dim_{\Fq} \Gamma_q(\Fq)_d = \! \sum_{j=2}^{m+1} (-1)^j {{m+1}\choose{j}} \sum_{i=0}^{j-2} {{d+(i+1)(q-1) - jq + m}\choose{d+(i+1)(q-1) - jq }}
\end{equation}
for any $d\ge 0$.
As noted in 
\S\,\ref{subsec:projred}, the space of projectively reduced polynomials in $\Fq[x_0,\dots,x_m]$ can be identified with $\Fq[x_0,\dots,x_m]/\Gamma_q(\Fq)$. In particular, $r_d \le {{m+d}\choose{d}}$ and the dimension of this space of projectively reduced polynomials in $\Fq[x_0,\dots,x_m]_d$
is $\binom{m+d}{d} - r_d$ for any nonnegative integer $d$. This also shows that
\begin{equation}\label{eq:Mbard-rd}
|\Monb_d| =  {{m+d}\choose{d}} - r_d \quad \text{for any $d\ge 0$.}
\end{equation}
Using this or otherwise (see, e.g.,  \cite[p. 237]{MR}), we readily see that
\begin{equation}\label{eq:rdvalues}
 r_d = 0 \text{ if } d\le q, \ r_{q+1} = {{m+1}\choose{2}}, \ \text{and} \
 r_d = {{m+d}\choose{d}} - \p_m \text{ if } d > m(q-1). 
 \end{equation}

\begin{definition}\label{def:ebar}
For any $d\ge 0$ and $1\le r \le \binom{m+d}{d} - r_d$,
we define
\begin{equation}\label{eq:ebardm}
\overline{e}_r(d,m):=\max_{G_1, \dots, G_r} |V(G_1, \dots, G_r)(\Fq)|,
\end{equation}
where the maximum is taken over all possible sets $\{G_1, \dots, G_r\}$ of $r$ linearly independent, projectively reduced polynomials in $\Fq[x_0,\dots,x_m]_d.$ As a natural convention, we set $\overline{e}_0(d,m):=0$.
\end{definition}

%
It is clear 
that $\overline{e}_r(d,m) \le e_r(d,m)$ for all $d\ge 0$ and
$1 \le r \le \binom{m+d}{d} - r_d$. 
A more precise relationship is given by the following.

 \begin{theorem}\label{ers}
Let $d, r$ be nonnegative integers such that $r \le \binom{m+d}{d} - r_d$. Then
$$
e_{r+r_d}(d,m)=\overline{e}_{r}(d,m).
$$
In particular, $e_r(d,m)=\overline{e}_r(d,m)$ if $d \le q.$
\end{theorem}

\begin{proof}
Fix a basis $\{\Phi_1, \dots , \Phi_{r_d}\}$ of the $\Fq$-vector space $\Gamma_q(\Fq)_d$. We will show that $e_{r+r_d}(d,m) \ge \overline{e}_{r}(d,m)$ and $e_{r+r_d}(d,m) \le \overline{e}_{r}(d,m)$.

Let $G_1, \dots, G_r \in \Fq[x_0, \dots, x_m]_d$ be linearly independent and projectively reduced polynomials such that $|V(G_1, \dots, G_r)(\Fq)|=\overline{e}_r(d,m).$ Then every nontrivial linear combination of $G_1, \dots, G_r$ is projectively reduced and therefore it does not belong to $\Gamma_q(\Fq)_d$. Consequently, the polynomials $G_1, \dots, G_r, \Phi_1, \dots, \Phi_{r_d}$ are linearly independent. Since $\Phi_1, \dots, \Phi_{r_d}$ vanish everywhere on $\PP^m (\Fq)$, we see that  $V(G_1, \dots, G_r)(\Fq)=V(G_1, \dots, G_r, \Phi_1, \dots, \Phi_{r_d})(\Fq)$.
It follows that
$$
e_{r+r_d}(d,m) \ge |V(G_1, \dots, G_r, \Phi_1, \dots, \Phi_{r_d})(\Fq)| = |V(G_1, \dots, G_r)(\Fq)| = \overline{e}_r(d,m).
$$

To prove the other inequality, suppose $F_1,\dots,F_{r+r_d}$ are linearly independent polynomials in
$\Fq[x_0,\dots,x_m]_d$ such that $|V(F_1,\dots,F_{r+r_d})(\Fq)|=e_{r+r_d}(d,m).$ Let $V_{r+r_d}$ denote the $\Fq$-linear span of $F_1,\dots,F_{r+r_d}$. Then
$$
\dim_{\Fq} \frac{V_{r+r_d} + \Gamma_q(\Fq)_d}{ \Gamma_q(\Fq)_d }
\ge \dim_{\Fq} V_{r+r_d} - \dim_{\Fq} \Gamma_q(\Fq)_d
= (r+ r_d) - r_d = r.
$$
So we can find $F'_1, \dots , F'_r \in V_{r+r_d}  + \Gamma_q(\Fq)_d$ that are  linearly independent (mod $\Gamma_q(\Fq)_d$). 
Let $G'_i = \overline{F}'_i$ denote the projective reduction of $F'_i$. Then $G'_1, \dots , G'_r$ are linearly independent and projectively reduced polynomials in $\Fq[x_0,\dots,x_m]_d$. Hence as in the previous paragraph, we see that $G'_1, \dots, G'_r, \Phi_1, \dots, \Phi_{r_d}$ are linearly independent.
 Also, $F'_i -G'_i \in  \Gamma_q(\Fq)$ for all $i=1, \dots , r$. 
It follows that
\begin{eqnarray*} 
e_{r+r_d}(d,m) & \ge & |V(G'_1, \dots, G'_r, \Phi_1, \dots, \Phi_{r_d} )(\Fq)| \\
&= & |V(G'_1,\dots,G'_r)(\Fq)|\\
&= & |V(F'_1,\dots,F'_r)(\Fq)|\\
 & \ge & |V(F_1,\dots,F_{r+r_d})(\Fq)| \\
 &=& e_{r+r_d}(d,m).
\end{eqnarray*}  
Consequently, equality holds throughout, and in particular,
$$
\overline{e}_r(d,m) \ge |V(G'_1,\dots,G'_r)(\Fq)|= e_{r+r_d}(d,m).
$$
This proves that $e_{r+r_d}(d,m)=\overline{e}_{r}(d,m)$. Finally, if $d\le q$, then $r_d=0$, and therefore $e_{r}(d,m)=\overline{e}_{r}(d,m)$ in this case.
\end{proof}

\begin{remark}
As noted in \cite[Rem. 6.2]{DG2}, one has $e_r(d,m)=\p_m$ if $0 \le r \le r_d.$ Indeed, it is clear that for any $r\ge 0$ one has $e_r(d,m) \le |\PP^m(\Fq)|=\p_m$, while for $r \le r_d$, one may choose linearly independent polynomials $G_1,\dots,G_r \in \Gamma_q(\Fq)_d$ and deduce that $e_r(d,m) \ge |V(G_1,\dots,G_r)(\Fq)|=\p_m.$
\end{remark}

The vanishing ideal of $\Aff^m(\Fq)$ is easy to determine; it is precisely the ideal of $\Fq[x_1, \dots , x_m]$ generated by $\{x_i^q - x_i : 1\le i \le m\}$. It is not difficult to see (using, e.g., the principle of inclusion-exclusion) that the dimension of the $\Fq$-vector space of polynomials of degree $\le d$ in this vanishing ideal is given by
$$
\rho_d:= \sum_{j=1}^m (-1)^{j-1} \binom{m}{j} \binom{m+d - jq}{d-jq} \quad \text{for all } d\ge 0.
$$
In case $d < q$, one can directly see that $\rho_d=0$.
Since the dimension of the $\Fq$-vector space of polynomials of degree $\le d$ in $\Fq[x_1, \dots , x_m]$ is clearly $\binom{m+d}{d}$, we see that $\rho_d \le \binom{m+d}{d}$ for all $d\ge 0$.

We have seen in \S\,\ref{subsec:Hr1} an affine analogue of $e_{r}(d,m)$, namely,
$e^{\Aff}_r(d,m):$. Here is a natural affine analogue of $\overline{e}_{r}(d,m)$. 
For $0\le r \le \binom{m+d}{d} - \rho_d$, we define
$$
\overline{e}^{\Aff}_r(d,m):=\max_{g_1, \dots, g_r} |Z(g_1, \dots, g_r)(\Fq)|,
$$
where the maximum is taken over $r$ linearly independent and reduced polynomials $g_1, \dots , g_r$ of degree $\le d$ in $\Fq[x_1, \dots , x_m]$. 
Arguing similarly as in the proof of Theorem~\ref{ers}, we see that
$$
e^{\Aff}_{r+\rho_d} (d,m) = \overline{e}^{\Aff}_r(d,m) \quad \text{for } 1\le r \le
\binom{m+d}{d} - \rho_d.
$$
Moreover, $e^{\Aff}_r(d,m) = \overline{e}^{\Aff}_r(d,m) $ if $d< q$, whereas $e^{\Aff}_r(d,m) = q^m$ if $r\le \rho_d$.

The result of Heijnen and Pellikaan \cite{HP} that was alluded to in the introduction solves the problem of determining $\overline{e}^{\Aff}_r(d,m) $.
A special case of this (when $d<q$) was stated earlier (Theorem~\ref{Hr1}). Here is the general version.

\begin{theorem}[Heijnen-Pellikaan]
\label{Hr}
Let $d$ be a nonnegative integer. Then
$$
\overline{e}^{\Aff}_r(d,m)  = H_r (d, m) \quad \text{for } 1\le r \le
\binom{m+d}{d} - \rho_d.
$$
\end{theorem}

\subsection{An Upper Bound}

In this subsection we will use the combinatorial results in Sections \ref{combin} and \ref{combin2} to obtain an upper bound for $e_r(d,m)$ when $d < q$. As we have seen in the last subsection, when $d \le q$, the
quantities $e_{r}(d,m)$ and $\overline{e}_{r}(d,m)$ coincide. Thus, we will work with $\overline{e}_{r}(d,m)$, and we begin by relating it with the ``maximal footprint" 
defined as follows.

\begin{definition}
Given any nonnegative integers $d, r$ and $e$ with $1\le r\le |\Monb_d|$, i.e., $1\le r \le  \binom{m+d}{d} - r_d$,  define
$$
A_r(d,m;e):= \max\{ |\Delta_e(\s)| : \s \subseteq \Monb_d \text{ with } |\s| = r\}.
$$
\end{definition}

The relation between $A_r(d,m;e)$ and $\overline{e}_{r}(d,m)$ is given by the lemma below. 
In the remainder of this section, we consider the lexicographic order $\lex$ on the set $\Mon$ of all monomials in $x_0, \dots , x_m$. For $0\ne F\in \Fq[x_0, \dots , x_m]$, the largest monomial (w.r.t. $\lex$) appearing in $F$ (with a nonzero coefficient) will be denoted by $\lm(F)$.

\begin{lemma}\label{lem:distinctLM}
Let $d, r$ be nonnegative integers with $1\le r\le   \binom{m+d}{d} - r_d$. Then
$$
\overline{e}_{r}(d,m) = \max_{F_1, \dots , F_r} |V(F_1, \dots , F_r)(\Fq)|,
$$
where the maximum is over all possible
sets $\{F_1, \dots , F_r\}$ of linearly independent and projectively reduced polynomials in $\Fq[x_0, \dots , x_m]_d$ such that $\, \lm(F_1), \dots, \lm(F_r)$ are distinct. Consequently,
$$
\overline{e}_{r}(d,m)  \le A_r(d,m;e) \quad \text{for all } e \gg 0.
$$
\end{lemma}

\begin{proof}
If $G_1, \dots , G_r\in \Fq[x_0, \dots , x_m]_d$ are linearly independent and projectively reduced, then we can easily obtain $F_1, \dots , F_r\in \Fq[x_0, \dots , x_m]_d$ that are linearly independent and projectively reduced such that $\lm(F_1), \dots , \lm(F_r)$ are distinct and $V(F_1, \dots , F_r) = V(G_1, \dots , G_r)$. For example, we can obtain them recursively by taking $F_1:=G_1$ and
for $1< i  \le r$, taking 
$F_i := G_i - c_1F_1 -  \dots c_{i-1} F_{i-1}$, where $c_1, \dots , c_{i-1}\in \Fq$ are chosen in such a way that none among 
$\, \lm(F_1), \dots , \lm(F_{i-1})$ appear in $F_i$. 
This proves the first assertion. The second assertion follows from the first using the projective $\Fq$-footprint bound 
(Theorem \ref{thm:qfootprint}).
\end{proof}
%

The following result could be viewed as a projective analogue of Theorem~\ref{thm:Wei}, which in turn, arose from the works of Clements--Lindstr\"om \cite{CL}, Wei \cite{W} and Heijnen-Pellikaan \cite{HP, H}. 
For any  integers $d, r$ with $d\ge 0$ and $1\le r \le \binom{m+d}{d} - r_d$,  we denote by 
$\m_d(r)$  the set of first $r$ elements of $\Monb_d$ in descending lexicographic order, and for any $d\ge 0$, we set
$\m_d(r)$ to be the empty set if $r=0$.
\begin{theorem}\label{thm:upperbound}
Let $d, r$ be integers with $1\le d < q$ and $1\le r \le \binom{m+d}{d} - r_d$.
Then
\begin{equation}
\label{eq:ProjWei}
|\Delta_e (\s)| \le \big|\Delta_e (\m_d(r))\big|  \quad \text{for all } e \gg 0 \text{ and }\s \subseteq \Monb_d \text{ with } |\s|=r.
\end{equation}
Consequently, $A_r (d, m; e) = |\Delta_e (\m_d(r))|$ and $\overline{e}_{r}(d,m)  \le |\Delta_e (\m_d(r))|$
for 
\hbox{$e \gg 0$.}
\end{theorem}

\begin{proof}
We prove \eqref{eq:ProjWei} by 
induction on $m$. Suppose $m=1$ and \hbox{$\s =\! \{\mu_1,\dots,\mu_r\} \! \subseteq \! \Monb_d$.} We may assume, without loss of generality, that $\mu_1 \glex \cdots \glex \mu_r$. Since 
the set $ \Monb_d$ consists 
of $x_0^d \glex x_0^{d-1}x_1 \glex \cdots \glex x_1^d$, we obtain 
 $\deg_{x_0} \mu_r < \cdots < \deg_{x_0} \mu_1$.

Now observe that if $e \gg 0$, then $\Delta_e^{(0)}(\s)=\emptyset$ if $x_0^d \in \s$, whereas   $\Delta_e^{(0)}(\s)=\{x_0^e\}$ if $x_0^d \not\in \s$. 
Also, if $i:=\deg_{x_0} \mu_r$ and $e\ge i-1$, then $\Delta_e^{(1)}(\s)=\{x_0^{i-1}x_1^{e-i+1},\dots,x_1^e\}$, and so 
 $|\Delta_e^{(1)}(\s)|=\deg_{x_0} \mu_r.$
Consequently,  if $x_0^d \in \s$, then $\deg_{x_0} \mu_r \le d-r+1$ and
$$
|\Delta_e(\s)|=|\Delta_e^{(0)}(\s)|+|\Delta_e^{(1)}(\s)|=0+\deg_{x_0} \mu_r \le d-r+1 \quad \text{for all } e \gg 0,
$$
whereas if $x_0^d \not \in \s$, then $\deg_{x_0} \mu_r \le d-r$ and
$$
|\Delta_e(\s)|=|\Delta_e^{(0)}(\s)|+|\Delta_e^{(1)}(\s)|=1+\deg_{x_0} \mu_r \le d-r+1 \quad \text{for all } e \gg 0.
$$
On the other hand,  $\m_d(r)=\{x_0^d,x_0^{d-1}x_1,\dots,x_0^{d-r+1}x_1^{r-1}\}$ and a similar reasoning shows that $|\Delta_e(\m_d(r))|=d-r+1$ for all $e \gg 0$.

Next suppose $m>1$ and \eqref{eq:ProjWei}  holds for all values of $m$ smaller than the given one.
Consider any $\s \subseteq \Monb_d$ with $|\s| = r$ and let $r':= |\s^{\langle m-1 \rangle}|$.
By Theorem \ref{thm:projectivetoaffine}, 
$$
|\Delta_e(\s)|=\sum_{\ell=0}^m |\Delta_e^{(\ell)}(\s)|=\sum_{\ell=0}^m |\FP^{(\ell)}(\sigma^{(\ell)}(\s^{\langle\ell\rangle}))|.
$$
Moreover $\s^{\langle m \rangle}=\s$, since $d<q$. This implies that $|\s| = |\f^{(m)} (\s^{\langle m \rangle})| = r$ and $\f^{(m)} (\s^{\langle m \rangle}) \subset \Nom_{\le d}^{(m)}$. Hence by Theorem \ref{thm:affinecomb}, 
$$
|\Delta_e^{(m)}(\s)|=|\FP^{(m)}(\sigma^{(m)}(\s^{\langle m \rangle}))| \le |\FP^{(m)}(T)|,
$$
where $T=\L_d^{(m)}(r') \cup \m_{d-1}^{(m)}(r-r')$.
Now define
$$
\t := \L_d^{(m)}(r') \cup x_m\m_{d-1}(r-r'),
$$
where $x_m\m_{d-1}(r-r')=\{x_m \mu \, : \, \mu \in \m_{d-1}(r-r')\}.$ 
Since $d<q$, we readily see that 
$\t \subset \Monb_d$, $|\t|=r$, and $\f^{(m)}(\t)=T.$ 
Therefore by Theorem~\ref{thm:projectivetoaffine},
$$
|\Delta_e^{(m)}(\s)|  \le |\FP^{(m)}(\f^{(m)}(\t) )| =  |\Delta_e^{(m)}(\t)|.
$$
On the other hand, by the induction hypothesis, 
$$
\sum_{\ell=0}^{m-1} |\Delta_e^{(\ell)}(\s)| \le \sum_{\ell=0}^{m-1} |\Delta_e^{(\ell)}(\t)|.
$$
Consequently, from equation \eqref{eq:disjointshfp} and the definition of $\t$, it follows that
\begin{equation}\label{eq:rearrange}
|\Delta_e(\s)| = \sum_{\ell=0}^m  |\Delta_e^{(\ell)}(\s)|  \le  \sum_{\ell=0}^m |\Delta_e^{(\ell)}(\t)|
=   |\Delta_e(\L_d^{(m)}(r') \cup x_m\m_{d-1}(r-r'))|.
\end{equation}
With this in view, we can replace $\s$ by $\L_d^{(m)}(r') \cup x_m \m_{d-1}(r-r')$. Thus, to prove \eqref{eq:ProjWei}, it suffices to show that $|\Delta_e(\s)| \le |\Delta_e(\t)|$ for all $e \gg 0$, where we now take
$$
\s :=\L_d^{(m)}(r') \cup x_m \m_{d-1}(r-r')  \quad \text{and} \quad \t:=\m_d(r) .
$$
In view of Remark \ref{rem:sigmam}, 
$\f^{(m)}(\t)=\m_d^{(m)}(r)$. Also clearly,  $\t^{\langle m-1 \rangle} = \L_d^{(m)}(s')$ and $\t \setminus \t^{\langle m-1 \rangle}=x_m\m_{d-1}(r-s')$
 for some 
 $s' \ge 0$. 
Now we distinguish two cases.

\smallskip

\noindent
{\bf Case 1:} $r' < s'.$

Let $\alpha$ be the $(r'+1)^{\rm th}$ element of $\Nom_d^{(m)}$. Since $r'<s'$, we see that $\alpha \in \L_d^{(m)}(s')$ and $\alpha \not \in \L_d^{(m)}(r').$ In particular, $\alpha \in \t \setminus \s.$ We claim that there exists $\beta \in \s$ such that $\beta \lex \alpha.$
Suppose, if possible, the claim is false. Then $\alpha \lexeq \beta$ for every $\beta \in \s$. Since $\a \in \t$ and $\t = \m_d(r)$, this will imply that $\beta \in \t$ for all $\beta\in \s$. Thus, $\s \subset \t$. Further, the fact that $|\s| = |\t| =r$ implies that $\s = \t$.
But this is a contradiction since 
 $\alpha \in \t \setminus \s.$
Hence the claim is true. Note that if $\beta \in \s$ is such that $\beta \lex \alpha$, then $\alpha \not \in \L_d^{(m)}(r')$ implies that $\beta \not \in \L_d^{(m)}(r')$, and so $\beta \in x_m \m_{d-1}(r-r')$.
Now choose $\gamma = x_0^{i_0}\cdots x_{m}^{i_{m}}$ to be the largest (in lexicographical order) element of $\s$ such that $\gamma \lex \alpha.$ Then, as noted above, $\gamma  \in x_m \m_{d-1}(r-r')$.
Consider $\mu: =  \gamma x_{m-1}/x_m$. Since $\gamma  \in x_m \m_{d-1}(r-r')$ and $d< q$, we see that $\mu \in \Monb_d$. Moreover, $\mu \glex \gamma$ and in fact, $\mu$ is the immediate successor of $\gamma$ in $ \Monb_d$ in the lexicographic order. Hence $\gamma \lex \alpha$ implies that $\mu \lexeq \alpha$. In case $x_m\mid \mu$, then $\mu  \in x_m \m_{d-1}(r-r')$ because $\mu \glex \gamma$ and  $x_m \m_{d-1}(r-r')$ is upwards closed in $x_m \Monb_{d-1}$. But then $\mu \in \s$, and this contradicts the maximality of $\gamma$. Thus $x_m \nmid \mu$, i.e., $i_m=1$ and $\mu \not \in \s$.  In particular, since $d< q$, we see that $i_{m-1}+1 < q$ and $x_0^{i_0}\cdots x_{m-2}^{i_{m-2}} x_{m-1}^{i_{m-1}+i_{m}} = \mu  \not \in \s$. Hence by Definition \ref{def:phi}, 
$\mu = \phi(\gamma) \in \phi(\s)^{\langle m-1 \rangle}$, even though $\gamma \not \in \s^{\langle m-1 \rangle}$. On the other hand, if $\nu \in \s^{\langle m-1 \rangle}$, then clearly, $\phi(\nu) = \nu$, and
so $\nu \in \phi(\s)^{\langle m-1 \rangle}$.
 This shows that $r'':=|\phi(\s)^{\langle m-1 \rangle}|>|\s^{\langle m-1 \rangle}|=r'.$ Moreover,
 by Theorem \ref{phi},  $|\Delta_e(\s)| \le |\Delta_e(\phi(\s))|.$ Now, if $r'' < s'$, then we iterate the procedure by replacing $\s$ with $\phi(\s)$. Else we proceed to the next case.

\smallskip

\noindent
{\bf Case 2:} $r' \ge s'.$

In this case, for  $\ell=0,\dots,m-1$, 
$$
\Delta_e^{(\ell)}(\s)=\Delta_e^{(\ell)}(\L_d^{(m)}(r')) \subseteq \Delta_e^{(\ell)}(\L_d^{(m)}(s'))=\Delta_e^{(\ell)}(\t) \quad \text{for all $e \gg 0$}.
$$
Moreover, $\f^{(m)}(\t)=\m_d^{(m)}(r)$ and so by Theorem \ref{thm:Wei},
$$
|\Delta_e^{(m)}(\s)|=|\FP^{(m)}(\f^{(m)}(\s^{\langle m \rangle}))| \le |\FP^{(m)}(\f^{(m)}(\t))|=|\Delta_e^{(m)}(\t)|
\quad \text{for all $e \gg 0$}.
$$
Hence  from equation \eqref{eq:disjointshfp}, we obtain $|\Delta_e(\s)| \le |\Delta_e(\t)|$ for all $e \gg 0$, as desired.
\end{proof}

%
%

We are now ready to obtain an upper bound for $e_r(d,m)$ mentioned 
in the introduction. For any integers $d, r$ with $d\ge 0$ and $1\le r\le \binom{m+d}{d} - r_d$,
let us define
$$
K_r(d,m):=  \sum_{j=0}^{m-1}a_j \p_{m-1-j},
$$
where $(a_0, \dots , a_m)$ denotes the $r^{\rm th}$ element, in lexicographic order, of the set of $(m+1)$-tuple $(b_0, \dots , b_m)$ of nonnegative integers such that $b_0+\cdots + b_m =d$, or in other words, $a_0, \dots , a_m$ are unique nonnegative integers such that  $x_0^{a_0}\cdots x_m^{a_m}$ is  the smallest monomial in $\m_d(r)$ in lexicographic order.

\begin{theorem} 
\label{prop:upper}
Let $d,r$ be integers such that $1\le d<q$ and $1 \le r \le \binom{m+d}{d}.$ Then $A_r (d, m;e) =K_r(d,m)$ for all $e\gg 0$. Consequently,
$$
e_r(d,m) \le K_r(d,m).
$$
\end{theorem}

\begin{proof}
Let $x_0^{a_0}\cdots x_m^{a_m}$ be the smallest monomial in $\m_d(r)$ in lexicographic order.
In view of Theorem \ref{thm:upperbound}, it suffices to show that
$$
|\Delta_e(\m_d(r))|  =   \sum_{j=0}^{m-1}a_j \p_{m-1-j} \quad \text{for all $e \gg 0$}.
$$
We will prove this by induction on $m$.
The case $m=1$ is easy, since we have already noted in the proof of Theorem \ref{thm:upperbound} that $|\Delta_e(\m_d(r))|=d-r+1=(d-r+1)\p_0$ and  $x_0^{d-r+1}x_1^{r-1}$ is the smallest monomial in $\m_d(r)$ in lexicographic order if $m=1$.
%

Now suppose that $m>1$ and the result holds for all values of $m$ smaller than the given one.
In view of 
Remark \ref{rem:sigmam}, 
$\sigma^{(m)}(\m_d(r)) = \m_d^{(m)}(r)$ 
and therefore 
from equation \eqref{eq:ssl} and Theorem \ref{thm:projectivetoaffine}, we see that
$$
|\Delta_e^{(m)}(\m_d(r))|=|\FP^{(m)}(\m_d^{(m)}(r))|=\sum_{j=0}^{m-1} a_j q^{m-1-j} \quad \text{for all $e \gg 0$},
$$
where the last equality follows from \cite[Prop. 5.9]{HP} (see also \cite[Lem. 4.2]{BD}).
On the other hand, from equations \eqref{eq:disjointshfp} and \eqref{eq:ssl}, we  see 
that 
\begin{eqnarray*}
|\Delta_e(\m_d(r))| &=& |\Delta_e^{(m)}(\m_d(r))|+\sum_{j=0}^{m-1}|\Delta_e^{(j)}(\m_d(r)^{\langle j \rangle})|  \\
&= & |\Delta_e^{(m)}(\m_d(r))|+ |\Delta_e ( \m_d(r)^{\langle m-1 \rangle})| \quad \text{for all $e \ge 0$},
\end{eqnarray*}
Note that the smallest monomial in $\m_d(r)^{\langle m-1 \rangle}$ equals $x_0^{a_0} \cdots x_{m-2}^{a_{m-2}} x_{m-1}^{a_{m-1} + a_m}.$
Hence using the induction hypothesis, we obtain
$$
|\Delta_e(\m_d(r))| = \sum_{j=0}^{m-1}a_j q^{m-1-j}+\sum_{j=0}^{m-2} a_j \p_{m-2-j}=\sum_{j=0}^{m-1}a_j \p_{m-1-j} \quad \text{for all $e \gg 0$}.
$$
This completes the proof.
\end{proof}

As an application of Theorem \ref{prop:upper}, 
we show below that Conjecture \ref{conj:extended} holds (in the affirmative) for some 
``large" values of $r$.

\begin{lemma}\label{lem:equal}
Let $d$ be a positive integer with $d < q$ and let $r =\sum_{a=1}^i \binom{m+d-a}{d-1}$ for some 
positive integer $i \le m+1$. Then Conjecture \ref{conj:extended} holds 
and
$$
e_r(d,m)=\p_{m-i}.
$$
\end{lemma}

\begin{proof}
If $i=m+1$, then by equation \eqref{Pascal2}, $r = \binom{m+d}{d}$, and in this case it is clear that  $e_r(d,m)=0=\p_{m-i}$. Now suppose $i\le m$ so that $r < \binom{m+d}{d}$. We claim that
the $(r+1)^{\rm th}$ monomial of $\Monb_d$ in descending lexicographical order is given by $x_i^d.$ This is clear if $i=m$, and the general case follows by decreasing induction on $i$ if we note that for $i< m$, the monomials $\mu \in \Monb_d$ that satisfy $x_i^d \glexeq \mu \glex x_{i+1}^d$ are precisely the monomials of degree $d$ in the $m-i+1$ variables $x_i, \dots , x_m$ that are divisible by $x_i$, and the number of such monomials is clearly $\binom{m-i+d-1}{d-1}$, i.e., $\binom{m+d-(i+1)}{d-1}$.
From the above claim, it follows that the $r^{\rm th}$ monomial of $\Monb_d$ in descending lexicographical order is given by $x_{i-1}x_m^{d-1}.$ Hence 
Theorem \ref{prop:upper} implies that  $e_r(d,m) \le \p_{m-1-(i-1)}=\p_{m-i}$, while 
Theorem \ref{lem:lower} shows 
that
$$
e_r(d,m) \ge \p_{m-i-1}+H_0(d-1,m-i)=\p_{m-i-1}+q^{m-i}=\p_{m-i}.
$$
This proves the lemma.
\end{proof}

In the next section, we use a little trick from coding theory and 
results proved in this section to prove the validity of Conjecture \ref{conj:extended} for some more 
values of $r$. 

\section{Connection with Projective Reed-Muller Codes}
\label{sec:codes}

We begin by recalling some basics about 
linear codes and the notion of generalized Hamming weight that is relevant for us. We will then consider the projective Reed-Muller codes, and show that the determination of their generalized Hamming weights is intimately related to the problem considered in this paper.
As before, $m$ will denote a fixed positive integer. Moreover, $n,k$ are positive integers with $k\le n$.

\subsection{Generalized Hamming weights and projective Reed-Muller codes}
Recall that a \emph{$q$-ary $[n,k]$ (error-correcting linear)} \emph{code} 
is defined as a $k$-dimensional $\Fq$-linear subspace of $\Fq^n$. 
If $C$ is a $q$-ary $[n, k]$ code, then the parameters $n$ and $k$ are known, respectively, as the \emph{length} and the \emph{dimension} of $C$.
For any 
$D \subseteq \Fq^n$, we define the \emph{support} of $D$ to be the subset 
$$
\supp(D) = \{i \in \{1,\dots, n\} :  c_i \neq 0 \ \mathrm{for \ some} \ (c_1,\dots, c_n) \in D\} 
$$
and the \emph{support weight} of $D$, denoted $\w_H(D)$, to be the cardinality of $\supp (D)$.

Let $C$ be a $q$-ary $[n, k]$ code. For $1 \le r \le k$, the $r^{\rm th}$ \emph{generalized Hamming weight} (also known as the $r^{\rm th}$ \emph{higher weight}) of $C$ is defined to be
$$
d_r(C):=\min \{\w_H(D) :  D \text{ subspace of } C \text{ with } \dim D=r\}.
$$
Note that $d_1 (C)$ is the minimum distance of $C$.
The notion of generalized Hamming weights was 
introduced by Wei \cite{W} 
and he showed that for any $q$-ary $[n, k]$ code $C$,
\begin{equation}\label{mon}
1 \le d_1(C) < d_2(C) \dots <d_k (C) \le n.
\end{equation}
Also, it is clear that $d_k(C)=n$ if $C$ is \emph{nondegenerate}, i.e., $C$ is not contained in a coordinate hyperplane of $\Fq^n$.

We now recall the projective Reed-Muller codes, 
introduced by Lachaud \cite{La}.
Let $P_1,\dots,P_{\p_m}$ be some fixed 
representatives  in $\Fq^{m+1}$ 
for the $\p_m$ points of $\PP^m(\Fq)$, e.g., 
we could represent each point of $\PP^m(\Fq)$ by
an  $(m+1)$-tuple of elements of $\Fq$, not all zero, such that the 
last nonzero coordinate is $1$. Consider 
the linear map
$$
\ev: \Fq[x_0,\dots,x_m] \to \Fq^{\p_m} \quad \makebox{given \ by} \quad \ev(F)=(F(P_1), \dots, F(P_{\p_m})).
$$
The $d^{\rm th}$ order \emph{projective Reed-Muller code}, denoted by $\PRM_q(d,m)$, is defined as the image of the space of homogeneous polynomials of degree $d$ in $x_0,\dots,x_m$ with coefficients in $\Fq$, under the map $\ev$, i.e.,
$\PRM_q(d,m):=\ev(\Fq[x_0,\dots,x_m]_d).$
Note that $\PRM_q(d,m)$ is 
a nondegeneate linear code of length $\p_m$. 
A formula for the dimension of $\PRM_q(d,m)$ can be deduced from the observations in 
\S\,\ref{subsec:projred} and \S\,\ref{subsec:RedEqns}. Indeed, these observations show 
that 
the kernel of the  map $\ev$ restricted to $\Fq[x_0,\dots,x_m]_d$ is precisely $\Gamma_q(\Fq)_d$; 
consequently, 
%
in view of equations \eqref{eq:rd} and \eqref{eq:Mbard-rd}, and with $r_d$ as in \eqref{eq:rd}, we see that
$$
\PRM_q(d,m) \cong \frac{\Fq[x_0,\dots,x_m]_d}{\Gamma_q(\Fq)_d} \quad \text{and} \quad
\dim \PRM_q(d,m) = |\Monb_d| = \binom{m+d}{d} - r_d.
$$
Thus we see that the dimension, say $k_q(d,m)$, of $\PRM_q(d,m)$ is given by
$$
k_q(d,m)= \binom{m+d}{d} - 
\sum_{j=2}^{m+1} (-1)^j {{m+1}\choose{j}} \sum_{i=0}^{j-2} {{d+(i+1)(q-1) - jq + m}\choose{d+(i+1)(q-1) - jq }}.
$$
In view of \eqref{eq:rdvalues}, this formula simplifies to $\binom{m+d}{d}$ if $d\le q$ and to
$\p_m$ if $d\ge m(q-1)+1$. In fact, if $d\ge m(q-1)+1$, then 
 $\PRM_q(d,m) = \Fq^{\p_m}$. Thus projective Reed-Muller codes are of interest only when $d\le m(q-1)$. The 
 following result of S{\o}rensen \cite[Thm. 1]{So} gives an alternative formula for the dimension and an explicit formula for the 
 minimum distance of projective Reed-Muller codes of an arbitrary order. 

\begin{theorem}[S{\o}rensen] 
\label{thm:soe}
Suppose $1 \le d \le m(q-1)$. Then the projective Reed-Muller code $\PRM_q(d,m)$ is a nondegenerate linear code with
\begin{enumerate}
\item[{\rm (i)}]
$\displaystyle{\dim \PRM_q(d,m) = \sum_{\substack{t=1 \\ t \equiv d \pmod{q-1} }}^d \left( \sum_{j=0}^{m+1} (-1)^j \binom{m+1}{j}\binom{t-jq+m}{t-jq} \right),}$
\item[{\rm (ii)}]
$d_1(\PRM_q(d,m)) = (q-s)q^{m-t-1},$ where $s$ and $t$ are  unique integers such that 
$d-1 = t(q-1) + s$ and $\, 0 \le s < q-1.$
\end{enumerate}
\end{theorem}

\subsection{Connection with Homogeneous Equations over Finite Fields} 
It was noted in \cite[\S\,4.2]{DG} that if $d\le q$, then $d_r(\PRM_q(d,m))=\p_m-{e}_r(d,m)$ for 
each positive integer $r \le \binom{m+d}{d}$. We observe 
that there is a more general relation.
\begin{lemma}
\label{lem:drer}
Let $d,r$ be any integers with $d\ge 1$ 
and $1\le r \le \binom{m+d}{d} - r_d$. Then
\begin{equation}\label{dr}
d_r(\PRM_q(d,m))=\p_m-\overline{e}_r(d,m).
\end{equation}
\end{lemma}

\begin{proof}
Clearly, for any $r$-dimensional subspace $D$ of  $\PRM_q(d,m),$ there exist $r$ linearly independent, projectively reduced polynomials $F_1,\dots,F_r \in \Fq[x_0,\dots,x_m]_d$
such that $D$ is the $\Fq$-linear span of $\ev(F_1), \dots , \ev(F_r)$,
and moreover, 
$$
\w_H(D)=\p_m-|\{i \, : \, F_j(P_i)=0\ \makebox{for all $1 \le j \le r$}\}|=\p_m - |V(F_1,\dots,F_r)(\Fq)|.
$$
This yields the desired equality.
\end{proof}

\begin{corollary}\label{cor:ermono}
Let $d$ be an integer satisfying $1\le d\le q$.
Then
$$
e_1(d,m) > e_2(d,m) > \dots > e_{\binom{m+d}{d}}(d,m) = 0.
$$
\end{corollary}

\begin{proof}
By Theorem \ref{ers}, if $d\le q$, then $\overline{e}_r(d,m) = {e}_r(d,m)$ for $1\le r \le \binom{m+d}{d}$. So the desired inequalities follow from equation \eqref{mon} and Lemma \ref{lem:drer}. It is obvious from the definition that $e_{\binom{m+d}{d}}(d,m) = 0$, or alternatively, it follows from noting that $\PRM_q(d,m)$ is always a nondegenerate code.
\end{proof}

%

The following result generalizes Lemma \ref{lem:equal} as well as \cite[Thm. 4.7]{DG}, and extends the validity of
Conjecture \ref{conj:extended} 
for some  additional values of $r$.

\begin{theorem}\label{eq'}
Suppose $1\le d < q$ and $r=\sum_{a=1}^i \binom{m+d-a}{d-1} - t$ for some
positive integer $i \le m+1$. 
and nonnegative integer $t < d$. 
Then
$$
e_r(d,m)=\p_{m-i}+t \quad \mathrm{and \ consequently,} \quad
d_r(\PRM_q(d,m))=\p_m-\p_{m-i}-t.
$$
\end{theorem}

\begin{proof}
By Lemma \ref{lem:equal}, $e_{r+t}(d,m)=\p_{m-i}$.  In particular, the result holds if $t=0$.
Now assume that $1\le t < d$. 
As in the proof of Lemma \ref{lem:equal}, we 
observe that the $r^{\rm th}$ monomial of $\Monb_d$ in descending lexicographical order is given by $x_{i-1}x_{m-1}^{t}x_m^{d-1-t}.$ 
Hence Theorem \ref{prop:upper} implies
that $e_r(d,m) \le \p_{m-i}+t$, whereas
Corollary \ref{cor:ermono} implies that $e_r(d,m) \ge t + e_{r+t}(d,m)=t+\p_{m-i}.$ This yields the desired results. 
\end{proof}

\begin{example}
Suppose $m=2$ and $d=4$. Assume that $q > 4$. First, we know from \cite[Thm. 5.3]{BDG} that
equation \eqref{IC} is valid for all $r \le \binom{m+2}{2}$. Thus we know $e_r(d,m)$ for $1\le r\le 6$. Next, Lemma \ref{lem:equal} covers the values $r= 10, 14$ and $15= \binom{2+4}{4}$, since $10 =  \binom{5}{3}$, while $14 =  \binom{5}{3} + \binom{4}{3}$ and
$15 = \binom{5}{3} + \binom{4}{3} + \binom{3}{3}$. The remaining values are taken care of by Theorem \ref{eq'}.  Consequently, we can write down  the complete set of generalized Hamming weights of
$\PRM_q(4,2)$ for $q\ge 5$. A similar conclusion holds when $(m, d) = (2, 3)$ and $q\ge 4$.
\end{example}

Finally, we remark that  since $\overline{e}_r(d,m) = {e}_r(d,m)$ when  $d\le q$,
it is clear from 
\eqref{dr} that the results and conjectures about $e_r(d,m)$ in Section \ref{sec:er} can be easily reformulated in terms of $d_r(\PRM_q(d,m))$. As a sampling of one such result, we give below
a  reformulation of Theorem \ref{prop:upper} 
in the spirit of that in \cite[Thm. 5.10]{HP}.

\begin{proposition}\label{prop:boundgenHweight}
Let $r, d$ be positive integers with  
$d<q$ and $r \le \binom{m+d}{d}.$ Further, let 
$\mathcal{Q}_{d}^m := \{(\a_0, \dots, \a_m) \in Q^{m+1} :  \sum_{i=0}^m \a_i = (m+1)(q-1) - d \}$, where as in~\S\,\ref{subsec:3.1}, $Q:=\{0,1 , \dots, q-1\}$.   If $(\beta_0,\dots,\beta_m)$ is the $r^{\rm th}$ element of $\mathcal{Q}_{d}^m$ in ascending lexicographic order, then
$$d_r(\PRM_q(d,m)) \ge m+1+\sum_{j=0}^{m-1} \beta_j \p_{m-1-j}.$$
\end{proposition}

\begin{proof}
Let $x_0^{a_0}\cdots x_m^{a_m}$ be the smallest monomial in $\m_d(r)$ in lexicographic order. Then  Lemma \ref{lem:drer} and
Theorem \ref{prop:upper} imply that
\begin{align*}
d_r(\PRM_q(d,m)) \ge & \  \p_m - \sum_{j=0}^{m-1}a_j \p_{m-1-j}\\
                   = & \  m+1 + \sum_{j=0}^{m-1}(q-1)\p_{m-1-j}- \sum_{j=0}^{m-1}a_j \p_{m-1-j}\\
                   = & \ m+1 + \sum_{j=0}^{m-1}(q-1-a_j)\p_{m-1-j}.
\end{align*}
The result now follows by noting that the $(m+1)$-tuple $(q-1-a_0,\dots,q-1-a_m)$ is the $r^{\rm th}$ element of $\mathcal{Q}_{d}^m$ in ascending lexicographic order.
\end{proof}

\begin{remark}
\label{rem:RVV}
In a recent work, Ramkumar, Vajha and Vijay Kumar \cite{RVV} have 
determined all the generalized Hamming weights of what they call the ``binary projective Reed-Muller code". However, the code they consider is not $\PRM_2(d,m)$ as defined above (and studied by Lachaud \cite{La}, S{\o}rensen \cite{So}, and others), but, in fact, a puncturing of a subcode of $\PRM_2(d,m)$. Indeed, they consider evaluations at points of $\PP^m(\FF_2)$ (which, in this case, is just $\FF_2^{m+1} \setminus\{(0,\dots, 0)\}$) of polynomials in $\FF_2[x_0, \dots , x_m]_d$ that are reduced in the affine sense. The resulting code is degenerate, in general, and so they puncture it suitably so as to obtain a nondegenerate code, say $C_2(d,m)$ for $1\le d\le m$. 
The length of $\PRM_2(d,m)$ is $2^{m+1}-1$, while that of $C_2(d,m)$ is $2^{m+1}-\sum_{i=0}^{d-1} \binom{m}{i}$. Likewise, the dimension of $\PRM_2(d,m)$ is $\binom{m+2}{2}$, while that of $C_2(d,m)$ is $\binom{m+1}{2}$. Evidently, the generalized Hamming weight $d_r(C_2(d,m))$, for which a formula is given 
in \cite{RVV}, 
provides an upper bound for 
$d_r(\PRM_2(d,m))$ when $1\le r \le \binom{m+1}{2}$, but the equality does not hold, in general. 
\end{remark}

\section*{Acknowledgements}

Peter Beelen would like to thank IIT Bombay where parts of
this work were carried out when he was there in January 2017 as a Visiting Professor.
Sudhir Ghorpade would like to thank the Technical University of Denmark for short visits 
in June 2017 and July 2018 when some of this work was done.

Peter Beelen gratefully acknowledges  the support from The Danish Council of Scientific Research (DFF-FNU) for the project \emph{Correcting on a Curve}, Grant No. 8021-00030B.  
Mrinmoy Datta is grateful for the support received 
from 
The Danish Council for Independent Research (Grant No. DFF-6108-00362) and the Research Council of Norway (Project No. 280731). 
 Sudhir Ghorpade gratefully acknowledges the support from  Indo-Norwegian Research grant INT/NOR/RCN/ICT/P-03/2018 from the Dept. of Science and Technology, Govt. of India and the Research Council of Norway (DST-RCN), MATRICS grant MTR/2018/000369 from the Science and Engineering Research Board, Govt. of India, and Award grant 12IRAWD009 from the Industrial Research and Consultancy Centre (IRCC), IIT Bombay. 
 
We are also thankful to the referee for some helpful comments and to Rati Ludhani for her careful proofreading of a preliminary version of this article.

\end{document}